\documentclass[12pt]{amsart}
\usepackage{amscd}
\usepackage{amssymb}
\newcommand{\C}{\mathbb C}
\newcommand{\Q}{\mathbb Q}
\usepackage[curve]{xypic}
\newcommand{\comment}[1]{\marginpar{\sffamily{\noindent\tiny #1
   \par}\normalfont}}
\renewcommand{\comment}[1]{}
\newbox\mybox
\def\overtag#1#2#3{\setbox\mybox\hbox{$#1$}\hbox to
  0pt{\vbox to 0pt{\vglue-#3\vglue-\ht\mybox\hbox to \wd\mybox
      {\hss$\scriptstyle#2$\hss}\vss}\hss}\box\mybox}
\def\undertag#1#2#3{\setbox\mybox\hbox{$#1$}\hbox to 0pt{\vbox to
    0pt{\vglue#3\vglue\ht\mybox\hbox to \wd\mybox
      {\hss$\scriptstyle#2$\hss}\vss}\hss}\box\mybox}
\def\lefttag#1#2#3{\hbox to 0pt{\vbox to 0pt{\vss\hbox to
      0pt{\hss$\scriptstyle#2$\hskip#3}\vss}}#1}
\def\righttag#1#2#3{\hbox to 0pt{\vbox to 0pt{\vss\hbox to
      0pt{\hskip#3$\scriptstyle#2$\hss}\vss}}#1}

\def\Dot{\lower.2pc\hbox to 2.5pt{\hss$\bullet$\hss}}
\def\Circ{\lower.2pc\hbox to 2.5pt{\hss$\circ$\hss}}
\def\Vdots{\raise5pt\hbox{$\vdots$}}
\def\splicediag#1#2{\xymatrix@R=#1pt@C=#2pt@M=0pt@W=0pt@H=0pt}

\renewcommand\frame[2][3pt]{\hbox{$\vcenter{\hbox{\vrule\vbox 
{\hrule\kern#1\hbox{\kern#1$#2$\kern#1}\kern#1\hrule}\vrule}}$}}
\newcommand\lineto{\ar@{-}}
\newcommand\dashto{\ar@{--}}
\newcommand\dotto{\ar@{.}}
\newcommand{\bt}{\bullet}
\newtheorem{theorem}{Theorem}[section]
\newtheorem{theoremi}{Theorem}
\newtheorem*{theorem*}{Theorem}
\newtheorem{proposition}[theorem]{Proposition}
\newtheorem{lemma}[theorem]{Lemma}
\newtheorem{corollary}[theorem]{Corollary}

\newtheorem{conjecture}[theorem]{Conjecture}
\newtheorem*{conjecture*}{Conjecture}
\newtheorem*{BStheorem*}{Bhupal-Stipsicz Theorem}

\theoremstyle{definition}
\newtheorem{remark}[theorem]{Remark}
\newtheorem*{example*}{Example}
\newtheorem*{examples*}{Examples}

\begin{document}
\title[Rational homology disk smoothings]
{On rational homology disk smoothings of valency $4$ surface singularities}
\author{Jonathan Wahl}
\thanks{Research supported under NSA grant no.\ H98230-08-1-0036}
\address{Department of Mathematics\\The University of North
  Carolina\\Chapel Hill, NC 27599-3250} \email{jmwahl@email.unc.edu}
\begin{abstract} Thanks to the recent work \cite{SSW} and \cite{bs}, one has a complete list of resolution graphs of weighted homogeneous complex surface singularities admitting a rational homology disk (``$\mathbb Q$HD") smoothing, i.e., one with Milnor number $0$.  They fall into several classes, the most interesting of which are the $3$ classes whose resolution dual graph has central vertex with valency $4$.  We give a uniform ``quotient construction" of the $\mathbb Q$HD smoothings for those classes; it is an explicit $\mathbb Q$-Gorenstein smoothing, yielding a precise description of the Milnor fibre and its non-abelian fundamental group.  This had already been done for two of these classes in \cite{SSW}; what is new here is the construction of the third class, which is far more difficult.  In addition, we explain the existence of two different $\mathbb Q$HD smoothings for the first class.

We also prove a general formula for the dimension of a $\mathbb Q$HD smoothing component for a rational surface singularity.  A corollary is that for the valency $4$ cases, such a component has dimension $1$ and is smooth.  Another corollary is that ``most" $H$-shaped resolution graphs cannot be the graph of a singularity with a $\mathbb Q$HD smoothing.  This result, plus recent work of Bhupal-Stipsicz \cite{bs}, is evidence for a general
\begin{conjecture*} The only complex surface singularities admitting a $\mathbb Q$HD smoothing are the (known) weighted homogeneous examples.
\end{conjecture*}
  
\end{abstract}
\maketitle
    
\begin{center}{\textbf{Introduction}}
    \end{center}
    
Suppose $f:(\mathbb C^3,0)\rightarrow (\mathbb C,0)$ is an analytic map germ such that $(V,0)\equiv (f^{-1}(0),0)$ has an isolated singularity at the origin.  Its \emph{Milnor fibre} $M$ is a "nearby fibre" of $f$ intersected with a small ball about the origin: $M=f^{-1}(\delta)\cap B_{\epsilon}(0)$ ($0<|\delta|<<|\epsilon|$).  $M$ is a  manifold of dimension $4$, compact with boundary the \emph{link} $L$ of the singularity, the intersection of $X$ with a small sphere.   $L$ is a compact oriented $3$-manifold, which can be reconstructed from any resolution of the singularity as the boundary of a tubular neighborhood of the exceptional set.   $M$ is simply-connected, and has the homotopy type of a bouquet of $\mu$ $2$-spheres, where the Milnor number $\mu=\text{rk}\ H_2(M)$ can be computed as the colength of the Jacobian ideal of $f$.  We may think of $M$ as a ``Stein filling" of $L$.

It is natural  to consider smoothings of (germs of) arbitrary normal complex surface singularities.  That is, one has a flat surjective map $f:(\mathcal V,0)\rightarrow (\mathbb C,0)$, where $(\mathcal V,0)$ has an isolated $3$-dimensional singularity, and the fibre $(f^{-1}(0),0)\equiv (V,0)$ is a normal surface singularity.  $V$ has a link $L$ as before, and L\^{e} and Hamm showed how to define a Milnor fibre $M$---again, a Stein filling of $L$.  $M$ need no longer be simply connected, though $b_1(M)=0$ \cite{gst}.  In fact, one
 could have
 $\mu(M)=0$, i.e., $M$ is a rational homology disk (``$\mathbb Q$HD"), with $H_i(M,\mathbb Q)=0$ for $i>0$.

In the late 1970's and early '80's, the author studied this last possibility and produced many examples, both published (\cite{w3}, \cite{lw}) and unpublished (but see \cite{ds}, p. 505; \cite{js}, p.123; \cite{k}).  Such a $(V,0)$ must be a \emph{rational surface singularity} (implying in particular that the resolution dual graph $\Gamma$ is a tree), with discriminant (=$|H_1(L)|$)  a square, and  resolution invariant $K\cdot K$ an integer.   Among cyclic quotient singularities (whose links are lens spaces), $\mathbb Q$HD smoothings exist exactly for those of type $p^2/pq-1$ (\cite{w3},(5.9.1) and \cite{lw},(5.10)).  (In the language of Koll\'ar -Shepherd-Barron \cite{ksb}, these are particular examples of ``$\mathbb Q$-Gorenstein smoothings," i.e., smoothings which are quotients of smoothings of the so-called index $1$ cover, which is Gorenstein).  We also constructed $\mathbb Q$HD smoothings for three triply-infinite families (later called $\mathcal W, \mathcal N, \mathcal M$), each corresponding to a spherical triple $(3,3,3), (2,4,4)$, or $(2,3,6)$, and each growing from the log-canonical singularity with graph respectively:
$$
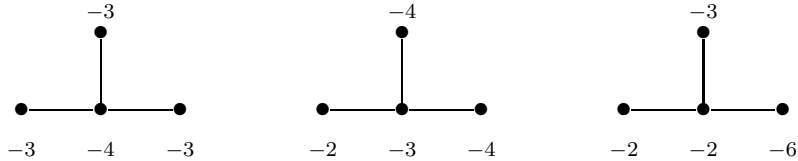
\xymatrix@R=6pt@C=24pt@M=0pt@W=0pt@H=0pt{
\\&&&&\\\
&\overtag{\bullet}{-3}{8pt}&&&&\overtag{\bullet}{-4}{8pt}&&&&\overtag{\bullet}{-3}{8pt}\\
&\lineto[u]&&&&\lineto[u]&&&&\lineto[u]\\
&\lineto[u]&&&&\lineto[u]&&&&\lineto[u]\\
&\lineto[u]&&&&\lineto[u]&&&&\lineto[u]\\
\undertag{\bullet}{-3}{4pt}\lineto[r]
&\undertag{\bullet}{-4}{4pt}\lineto[u]\lineto[r]&\undertag{\bullet}{-3}{4pt}
&&\undertag{\bullet}{-2}{4pt}\lineto[r]
&\undertag{\bullet}{-3}{4pt}\lineto[u]\lineto[r]&\undertag{\bullet}{-4}{4pt}
&&\undertag{\bullet}{-2}{4pt}\lineto[r]
&\undertag{\bullet}{-2}{4pt}\lineto[u]\lineto[r]&\undertag{\bullet}{-6}{4pt}\\
&&&&\\
&&&&\\
&&&&
}
$$
(The $\mathcal W$ family is \cite{w3}, 5.9.2).  We found other one-parameter infinite families; all were weighted homogeneous singularities, hence (\cite{ow}) $\Gamma$ is star-shaped (a \emph{star} is a vertex of valency at least $3$).  Most had a star of valency three; but there were also valency $4$ families for each triple.  We list the graphs below for each $p\geq 2$, and name the type via a triple $(a,b,c;d)$:
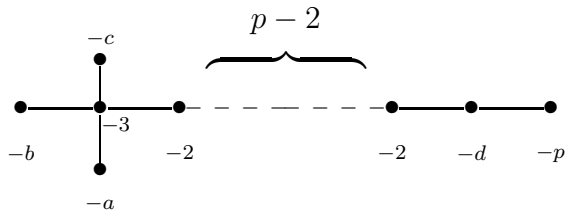
$$
\xymatrix@R=6pt@C=24pt@M=0pt@W=0pt@H=0pt{
\\&&&&p-2\\
&&\overtag{\bt}{-c}{8pt}
&&{\hbox to 0pt{\hss$\overbrace{\hbox to 60pt{}}$\hss}}&\\
&&\lineto[u]\\
&\undertag{\bt}{-b}{6pt}\lineto[r]&
\bt\lineto[u]\lineto[r]_(.2){-3}&\undertag{\bt}{-2}{6pt}\dashto[r]&\dashto[r]
&\undertag{\bt}{-2}{6pt}\lineto[r]
&\undertag{\bt}{-d}{6pt}\lineto[r]
&\undertag{\bt}{-p}{6pt}\\
&&\lineto[u]\\
&&\lineto[u]\\
&&\undertag{\bt}{-a}{2pt}\lineto[u]&&&\\
&&&&\\
&&&&\\
}
$$

\begin{align*}
\mathcal A^4\ &:\ (a,b,c;d)=(3,3,3;4)\\
\mathcal B^4\ &:\ (a,b,c;d)=(2,4,4;3)\\
\mathcal C^4\ &:\ (a,b,c;d)=(2,3,6;2)
\end{align*}
The main goal of the current paper is to give a new and uniform ``quotient construction" of a $\mathbb Q$HD smoothing of these valency $4$ examples, allowing an explicit description of the $\mathbb Q$HD Milnor fibre and its fundamental group.  (The $\mathcal A^4$ and $\mathcal B^4$ examples are easier and were already written down in \cite{SSW}, (8.12) and (8.14).) 

There are two ways to construct smoothings.  The first is H. Pinkham's method of ``smoothing with negative weight"  for a weighted homogeneous singularity \cite{p3}.  Here, the existence of a smooth projective variety $\tilde{X}$ with a specific curve configuration can be used to produce a smoothing of $V$, as proved carefully in \cite{SSW}, Theorem 8.  (One can also use this method to prove non-existence of $\mathbb Q$HD smoothings in certain cases, cf.\cite{bs}.)  But the total space of the smoothing has coordinate ring $\bigoplus \Gamma (\tilde{X},\mathcal O(D_k))$, for some divisors $D_k$ on $\tilde{X}$, so it can be difficult to deduce, for instance, the fundamental group of the Milnor fibre, or whether the smoothing is $\mathbb Q$-Gorenstein.  

One would prefer a ``quotient" construction (\cite{w3}, 5.9), which would be more explicit and give extra information.  In this case, one starts with (a germ of) an isolated $3$-diimensional Gorenstein singularity $(\mathcal Z,0)$; a finite group $G$ acting on it, freely off the origin; and a $G$-invariant function $f$ on $\mathcal Z$, whose zero locus $(W,0)$ has an isolated singularity (hence is normal and Gorenstein).  Then  $f:(\mathcal Z,0)\rightarrow (\mathbb C,0)$ is a smoothing of $W$, with Milnor fibre $M$; in our examples, $M$ will be simply-connected.  But one also has  $f:(\mathcal Z/G,0)\rightarrow (\mathbb C,0)$,  giving a $\mathbb Q$-Gorenstein smoothing of $(W/G,0)\equiv(V,0)$, with Milnor fibre $M/G$ (note $G$ acts freely on $M$).   If  the order of $G$ is equal to the Euler characteristic of $M$, then $M/G$ will have Euler characteristic $1$, hence will be a $\mathbb Q$HD Milnor fibre.  The resolution dual graph of $W/G$ can be calculated by resolving $\mathcal Z$ or $W$ and dividing by the induced action of $G$. 

A quotient construction is given in \cite{SSW}, (8.2) for the cyclic quotients and for types $\mathcal W$ and $\mathcal N$; in these cases, $G$ is abelian.  In addition, the  valency $4$ examples of types $\mathcal A^4$ and $\mathcal B^4$ are constructed there (modulo typos in (8.14): $m$ should be $2p^2-2p+1$, and the graph should have $p$ replaced by $p-2$).  In these cases, $(\mathcal Z,0)$ is respectively $(\mathbb C^3,0)$ or a quadric $(\{X_1X_3-X_2X_4=0\}\subset \mathbb C^4,0)$; $G$ is a metacyclic group; $f$ is an equation of ``Klein type"; and $M$ is simply connected, so that the fundamental group of the $\mathbb Q$HD smoothing is the non-abelian group $G$.  The new result in (7.4) extends this construction to the case $\mathcal C^4$:

 \begin{theoremi} For $p\geq 2$, write $m=p^2-p+1$, $n=6p$, $r=m-p+1$, and $\zeta$ = exp$(2\pi i/m)$, $\eta$=exp$(2\pi i/n)$, $\gamma= \text{exp}(2\pi i/6)$.  Consider
 \begin{enumerate}
 \item $\mathcal Z \subset \mathbb C^7$  the cone over an appropriate del Pezzo surface of degree $6$ in $\mathbb P^6$
 \item $G\subset \text{GL}(7,\mathbb C)$ the metacyclic group generated by 
 $$S(a_1,a_2,\cdots,a_7)=(\zeta a_1,\zeta^{r} a_2,\cdots,\zeta^{r^{5}}a_6,a_7)$$ 
 $$T(a_1,a_2,\cdots,a_7)= (\eta a_2,\cdots, \eta a_{6},\eta a_1,\eta a_7).$$
 \item $f=X_1^{p-1}X_2+\gamma X_2^{p-1}X_3+\gamma^2X_3^{p-1}X_4+\gamma^3X_4^{p-1}X_5+\gamma^4X_5^{p-1}X_6+\gamma^5X_6^{p-1}X_1.$ 
 \end{enumerate} Then the quotient construction applied to this data yields  $\mathbb Q$HD smoothings of type $\mathcal C^4$, whose Milnor fibre has fundamental group $G$.

 \end{theoremi}
 
   The total space for each of these valency $4$ smoothings is the quotient of the cone $\mathcal Z$ of an embedding $Z\subset \mathbb P^n$ for which $-K_Z$ is cut out by a hypersurface, and $f\in\Gamma(Z,-pK_Z)$.   (These three $Z$ are precisely the toric del Pezzo surfaces with cyclic symmetry, of orders $3,4,6.$)   $\mathcal Z$ has a canonical Gorenstein singularity, so the total spaces of the smoothings are \emph{log-terminal}; in a forthcoming paper, we prove this is a general phenomenon.
   
The relevance of $\mathbb Q$HD smoothings to symplectic topology arose from a 1997 paper of Fintushel-Stern  \cite{fs},  which showed that  one could ``blow-down" on a symplectic four-manifold a configuration of $2$-spheres that numerically correspond to a resolution of a cyclic quotient singularity of type $p^2/(p-1).$  Roughly speaking, one removes an appropriate symplectic neighborhood of the configuration, with boundary the lens space, and then ``pastes in the $\mathbb Q$HD Milnor fibre".  (This interpretation came later.)  A key point is that Seiberg-Witten invariants can be controlled, so one can create new smooth $4$-manifolds.   Another important consequence was the construction by Y. Lee and J. Park \cite{lp} of new minimal surfaces of general type, by globalizing $\mathbb Q$HD smoothings of the $p^2/pq-1$ cyclic quotient singularities.

The point of view of ``nice" filling of compact $3$-manifolds by $\mathbb Q$HD Milnor fibres was due to the symplectic topologists A. Stipsicz and Z. Szab\'o, leading to our collaboration  \cite{SSW} in 2007.  
The main achievement in \cite{SSW} was to restrict greatly the possible resolution diagram of a $(V,0)$ which could admit a $\mathbb Q$HD smoothing, and to list many examples where such smoothings exist.  (The same restrictions applied to similar problems in ``symplectic filling.")  The possible diagrams were of $7$ types.  The first class of graphs ($\mathcal G$) correspond to the lens spaces $p^2/(pq-1)$.  The next three ($\mathcal W, \mathcal N, \mathcal M$) are the 3 triply-infinite families mentioned above; in particular, all graphs in the first $4$ classes can be realized by singularities with $\mathbb Q$HD smoothings.  But the final 3 classes ($\mathcal A, \mathcal B, \mathcal C$) contain some diagrams of singularities that could not have a $\mathbb Q$HD smoothing (e.g., the singularity might not be rational).  Nonetheless, these classes contained six one-parameter valency $3$ star-shaped families, as well as the three one-parameter families $\mathcal A^4, \mathcal B^4, \mathcal C^4$, all of which correspond to singularities having a $\Q$HD smoothing.  (These families were known to the authors of  \cite{SSW},  but not all were listed there.)    An important consequence of \cite{SSW} was to give a neat way to organize all the families on our old list.  This result culminated in recent work of Bhupal-Stipsicz \cite{bs}, which says that for star-shaped graphs the old list is complete:

\begin{BStheorem*}[\cite{bs}] The examples described above are the only star-shaped graphs of singularities with a $\mathbb Q$HD smoothing. 

\end{BStheorem*}

Thus, \cite{bs} considered all the other star-shaped graphs allowed by Theorem 1.8 in \cite{SSW} (necessarily of classes $\mathcal A, \mathcal B$ or $ \mathcal C$), and proved (using ideas related to smoothing of negative weight) that a singularity with such a graph could not have a $\mathbb Q$HD smoothing.  In fact, \cite{bs} worked in the symplectic category and proved a stronger result, that the relevant $3$-manifolds did not even have a ``$\mathbb Q$HD weak symplectic filling'' .

The Bhupal-Stipsicz result lists the graphs, so one should say as much as possible about the singularities themselves and their $\mathbb Q$HD smoothings.  The valency $3$ examples which occur are all \emph{taut} in the sense of H. Laufer \cite{l}: there is a unique singularity with that graph,  so it is the  weighted homogeneous one (\cite{p2}).  But one should decide how many $\mathbb Q$HD smoothings it has.   For the valency $4$ graphs, \cite{l} proves that for each cross-ratio of $4$ points on the central curve, there is a unique analytic type (necessarily weighted homogeneous).  The quotient $\mathbb Q$HD smoothings in each of the cases of $\mathcal A^4, \mathcal B^4, \mathcal C^4$ work for only one particular cross-ratio, which in the first two cases can be identified because of the existence of a special symmetry of the $4$ points.  A general result (Theorem 8.1) gives a formula for the dimension of a $\mathbb Q$HD smoothing component.  As a Corollary, one has (see (8.2) below) the

\begin{theoremi} For a weighted homogeneous surface singularity, a $\mathbb Q$HD smoothing component in the base space of the semi-universal deformation has dimension one.
\end{theoremi}

It follows that there are at most finitely many cross-ratios for singularities with a $\mathbb Q$HD smoothing, and perhaps the calculations of \cite{bs} will imply these are the only ones.  
 Further, in all examples made from the quotient construction one finds such a $\mathbb Q$HD smoothing component is smooth.
 
 On the other hand, we prove (Theorem 2.4) that for type $\mathcal A^4$,  there are \emph{two}  $\mathbb Q$HD smoothing components in the base space of the semi-universal deformation of the singularity (this was also hinted at by J. Stevens in \cite{js}, p.124).   From one point of view, the two different smoothings come from two $3$-dimensional fixed-point free representations of the same group, which do not differ by an automorphism of the group (cf. Theorem \ref{A} below).   These two smoothings have non-isomorphic three-dimensional total spaces, but the two Milnor fibres are apparently diffeomorphic.

It is natural to formulate the following

\begin{conjecture*} Suppose a rational surface singularity has a $\mathbb Q$HD smoothing.  Then the singularity is weighted homogeneous, hence is one of the known examples.
\end{conjecture*}

To prove this Conjecture, by \cite{SSW} one would need to rule out those graphs of types $\mathcal A, \mathcal B, \mathcal C$ with two or more stars (each of which must have valency exactly $3$, except perhaps for an ``initial" one of valency $4$).  We  give a simple proof of non-existence of $\mathbb Q$HD smoothings for most of the candidate graphs with two valency $3$ stars.

\begin{theoremi}  Graphs of type $\mathcal A, \mathcal B, \mathcal C$ with two valency $3$ stars have resolution dual graph
$$
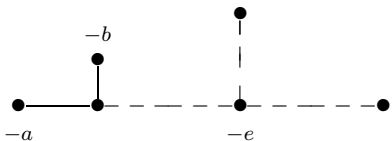
\xymatrix@R=12pt@C=24pt@M=0pt@W=0pt@H=0pt{
&&&&{\bullet}\dashto[dd]\\
&&\overtag{\bullet}{-b}{10pt}&&\dashto[u]\\
&\undertag{\bullet}{-a}{0pt}\lineto[r]&
{\bullet}\lineto[u]\dashto[r]&\dashto[r]&
\undertag{\bullet}{-e}{0pt}\dashto[r]&\dashto[r]&{\bullet}\\
&&&&\\
}
$$
with $e\geq 2$.  If $e\geq 3$, then any singularity with this graph cannot admit a $\mathbb Q$HD smoothing.
\end{theoremi}
 
 In the first section, we recall from J. Wolf's book \cite{wolf} the metacyclic groups $G$ which have faithful irreducible representations in $\mathbb C^d$ acting freely off the origin.  We also explain how two different such representations of the same group could yield non-isomorphic quotient singularities $\mathbb C^d/G$.  Sections $2$ and $3$ discuss cases $\mathcal A^4$ (where $d=3$) and $\mathcal B^4$ (where $d=4$), respectively; one studies conjugacy classes in $G$ and isotropy for the projective representation, to describe the singularities made via appropriate spaces $\mathcal Z$ and Kleinian polynomials $f$.  The difficult case $\mathcal C$ requires in Section $4$ a description of an appropriate
 $G$ (for $d=6$) but now with a \emph{reducible} action on $\mathbb C^7$;  in Section $5$ a del Pezzo $Z$ in $\mathbb P^6$ on whose cone $G$ acts, and whose non-trivial isotropy is described;  in Section 6 the $G$-invariant Kleinian polynomial cutting out a smooth curve $D$ on $Z$; and finally in Section $7$ a description of the smoothing and properties of its $\mathbb Q$HD Milnor fibre.  The last Section 8 discusses $\mathbb Q$HD smoothing components, and makes the aforementioned conjecture about the non-existence of other examples; this section is independent of all the others.
 
 We have profited from  conversations with Eduard Looijenga, Walter Neumann, and Andr\'{a}s Stipsicz.

\section{Some fixed-point free group actions}
\bigskip
Section 5.5 in J. Wolf's book \cite{wolf} discusses metacyclic groups acting freely off the origin of $\C^d$.  Consider positive integers $m,n,r, d, n'$, with $m>1$ and
\begin{enumerate}
\item $n=n'd$ 
\item $((r-1)n,m)=1$
\item $r$ is a multiplicative unit in $\mathbb Z/(m)$, of order $d$ 
\item every prime divisor of $d$ divides $n'$.  
\end{enumerate}
(In the cases of interest to us, $d=3,4,\text{or} \ 6$.)  Let $G$ be the group defined by generators $A$ and $B$ satisfying $$A^m=B^n=1,\ \\ BAB^{-1}=A^r.$$  $G$ is the semi-direct product of the cyclic groups $<A>$ and $<B>$, hence has order $mn$.  Its abelianization is cyclic of order $n$, generated by the image of $B$.

Suppose $k$ and $l$ are integers satisfying $(k,m)=(l,n)=1$.  Define  $d$-dimensional (unitary) representations $\pi_{k,l}$ of $G$ as follows:
Let  $\zeta$ = exp$(2\pi i/m)$ and $\eta$ = exp$(2\pi i/n)$.   Set
  $$\pi_{k,l}(A)(a_1,a_2,\cdots,a_d)=(\zeta^k a_1,\zeta^{kr}a_2,\cdots,\zeta^{kr^{d-1}}a_d)$$ 
 $$\pi_{k,l}(B)(a_1,a_2,\cdots,a_d)=(a_2,\cdots, a_{d},\eta^{ld}a_1).$$
 These give representations $\pi_{k,l}:G\rightarrow GL(d,\C)$, which depend only on $k$ mod $m$ and $l$ mod $n'$. 
 \begin{theorem} \cite{wolf}
 \begin{enumerate}
 \item The $\pi_{k,l}$ are irreducible faithful complex representations of $G$,  and act freely on $\mathbb C^d$ off the origin
 \item The subgroups $\pi_{k,l}(G)$ and $\pi_{1,l}(G)$ are equal
 \item $\pi_{1,l}$ is equivalent to $\pi_{1,l'}$ iff $\pi_{1,l}=\pi_{1,l'}$ iff $l'\equiv l$ modulo $n'$
 \item $\pi_{1,l}$ is equivalent to $\pi_{1,l'}\circ \varphi$, with $\varphi$  an automorphism of $G$, iff there exists $t\equiv1$ mod $d$ such that $\pi_{1,l'}=\pi_{1,tl}$
 \end{enumerate} 
 \begin{proof}  The statements are easy or taken directly from Theorem 5.5.6 of \cite{wolf}.   
 \end{proof}
 \end{theorem}
  We clarify a count not expressly made in the Remark on p. 170 of \cite{wolf}.
 \begin{corollary} $\pi_{1,l}$ differs from $\pi_{1,l'}$ by an automorphism of $G$ iff $l'\equiv l$ mod $GCD(d,n')$.
 \begin{proof} By the Theorem, one needs to examine the conditions that there is a $t\equiv 1$ mod $d$ such that $l'\equiv tl$ mod $n'$.  Reducing the last equivalence mod $GCD(d,n')$ gives that $l'\equiv tl\equiv l$ mod this $GCD$.  Conversely, suppose that $l,l'$ are prime to $n'$ and are congruent to each other mod $GCD(d,n')$.  Then there is an integer $s$ (prime to $n'$) so that $l'\equiv sl$ modulo $n'$, and $s\equiv1$ mod $GCD(d,n')$.  It is an elementary exercise to show one can find an $x$ so that $t=s+n'x$ is congruent to $1$ modulo $d$.
 \end{proof}
 \end{corollary}
 Next,  the images of $G$ under two representations are conjugate subgroups if and only if the corresponding representations become equivalent via an automorphism of $G$. 
Let  $G_l \subset GL(d,\C)$ be the image $\pi_{1,l}(G).$  
\begin{corollary}
\begin{enumerate}
\item  The number of non-conjugate subgroups $G_l$ is $\phi(GCD(n',d))$.
\item  The number of analytically distinct germs $(\C^d/G_l,0)$ is $\phi(GCD(n',d))$.
\end{enumerate}
\begin{proof}  The first statement is an immediate consequence of the last Corollary. For the second, a local isomorphism lifts to a local analytic isomorphism of $\C^d-\{0\}$ which is equivariant with respect to the actions of the two subgroups.  The automorphism extends over the origin, and its linearization intertwines the two group actions.
\end{proof}
\end{corollary}

\begin{remark} Note $GCD(n',d)>1$, so $\phi (GCD(n',d))$ is even unless the $GCD$ equals $2$.  In the latter case, there is only one subgroup $G_l$.  In all other cases, $G_1$ and $G_{-1}$ are non-conjugate.   Although the quotient singularities $\C^d/G_l$ are not analytically equivalent, the quotients of the spheres $S^{2d-1}/G_1$ and $S^{2d-1}/G_{-1}$ are diffeomorphic, in a strong sense.  For this question, one has to look at the orthogonal representations $\hat{\pi}_{k,l}:G\rightarrow O(d)$ as described in Theorem 5.5.11 of \cite{wolf}.  We summarize by saying that $S^{2d-1}/G_l$ and $S^{2d-1}/G_{l'}$ are isometrically equivalent iff $l\equiv\pm l'$ modulo $GCD(n',d)$.
\end{remark}

From now on, we use slightly different transformations to describe the $\pi_{k,l}$.  The diagonal transformation $$R(a_1,a_2,\cdots,a_d)=(a_1,\eta^{-l}a_2,\cdots,\eta^{-l(d-1)}a_d)$$ 
 commutes with $\pi_{k,l}(A)$, while
$$R^l\pi_{k,l}(B)R^{-l}(a_1,a_2,\cdots,a_d)=(\eta^la_2,\cdots,\eta^la_d,\eta^la_1).$$
Thus,  $\pi_{k,l}$  is equivalent (via $R^l$) to the representation defined by the transformations $S_k$ and $T_l$, where
\begin{equation*}
\tag{1.4.1}\label{1}S_k(a_1,a_2,\cdots,a_d)=(\zeta^{k}a_1,\zeta^{kr}a_2,\cdots,\zeta^{kr^{d-1}}a_d)
\end{equation*}
\begin{equation*}
\tag{1.4.2}\label{2} T_l(a_1,a_2,\cdots,a_d)= (\eta^{l}a_2,\cdots, \eta^{l}a_{d},\eta^{l}a_1).
 \end{equation*}
We will frequently write $S=S_1,\ T=T_1$. 

 Since the transformations are unitary, the transpose inverse is the complex conjugate; thus, the dual representations satisfy
 $$\pi_{k,l}^*\cong\pi_{-k,-l}.$$ 
 In particular, as $T(e_i)=\eta e_{i-1}$ (where the $e_i$ are standard basis vectors and $e_0\equiv e_d$), we have as corresponding action on the coordinate functions $T(X_i)=\eta^{-1}X_{i-1}$  (where $X_0\equiv X_d$).

\begin{lemma}  Let $G$ act on $\mathbb C^d$ via $\pi_{k,l}$, as in \eqref{1} and \eqref{2} above.  Suppose that $n'=m-r+1$, and write $\gamma=\eta^{l(m-r+1)}$.  Then 
$\gamma$ is a primitive $d^{th}$ root of 1, and 
 the polynomial $$f(X_1,X_2,\cdots,X_d)=X_1^{m-r}X_2+\gamma X_2^{m-r}X_3+\cdots +\gamma^{d-1}X_d^{m-r}X_1$$
is $G$-invariant.
\begin{proof} A straightforward calculation.
\end{proof}
\end{lemma}

We note a simple calculation from the defining relations of $G$: for $k\geq 0,$ all $j$,
 $$A^{-j}B^kA^j=A^{j(r^k-1)}B^k.$$
 
 \begin{lemma} If $i$ is a multiple of $r^k-1$ in $\mathbb Z/(m)$, then $A^iB^k$ is conjugate to  $B^k$.
 \begin{proof}  The hypotheses imply there is a $j\geq 0$ so that $i\equiv j(r^k-1)$ modulo $m$.  Now use the calculation.
 \end{proof}
 \end{lemma}
 
 \section{Class $\mathcal A$ valency 4 smoothings}
\bigskip

We consider the set-up of the last section for some three-dimensional representations, i.e. $d=3$.   Fix an integer $p\geq 2$, and write $m=3p^{2}-3p+1,  \ n=9p, \ n'=3p$, and $r=m-(3p-1)$. 
  Let  $\zeta$ = exp$(2\pi i/m)$,  and $\eta$ either exp$(2\pi i/n)$ or its inverse.  Consider the subgroup
 group $G_\eta =G\subset GL(3,\C)$ generated by the transformations $S$ and $T_\eta=T$, where 
 $$S(a,b,c)=(\zeta a,\zeta^{-(3p-1)}b,\zeta^{3p-2}c)$$ 
 $$T(a,b,c)= (\eta b, \eta c,\eta a).$$
In the notation of the first section, this subgroup arises  from the representation $\pi_{1,1}$ or $\pi_{1,-1}$, depending on which root of unity is $\eta$. 
%
 
 It is straightforward to verify the following:
 \begin{enumerate}
 \item $T^3$ is scalar multiplication by $\eta^{3}$
 \item The image $\bar{G}$ of $G$ in PGL$(3,\mathbb C)$ is $G/<T^3>$, hence has order $3m$ and is generated by $\bar{S}$ of order $m$ and $\bar{T}$ of order $3$.
 \end{enumerate}
 Checking that $r-1$ and $r^2-1$ are units in $\mathbb Z/(m)$,  Lemma 1.6 yields
\begin{lemma}
If $k$ is $1$ or  $2$,  then $S^iT^k$ is conjugate to $T^k.$
\end{lemma}
 
 \begin{lemma}
\begin{enumerate}
\item $\bar{G}$ acts faithfully on $\mathbb P^2$, with non-trivial isotropy only at $[1,0,0]$, $[0,1,0]$, $[0,0,1]$ and the $\bar{G}$-orbits of  the three points $[1,\rho,\rho^2]$, where $\rho$ is any cube root of $1$.
\item $<\bar{S}>$ is the $\bar{G}$-stabilizer of each of the $3$ points  $[1,0,0], [0,1,0], [0,0,1]$, and $<\bar{T}>$ cyclically permutes them.
\item $<\bar{T}>$ is the $\bar{G}$-stabilizer of each of the $3$ points  $[1,\rho,\rho^2]$, and these points are in distinct $\bar{G}$-orbits.
\end{enumerate}
\begin{proof}
First check the action of $\bar{G}$ on $\mathbb P^2$.  By the last Lemma, it suffices to find points fixed by $\bar{S}^i$, $\bar{T}$, and $\bar{T}^2=\bar{T}^{-1}$ (so the last case is unnecessary).  The calculations concerning $\bar{G}$ are now straightforward. 
\end{proof}
\end{lemma}

 Next, $G$ acts on the 3 coordinate functions on $\C^3$ 
 by $$S(X_1,X_2,X_3)=(\zeta^{-1} X_1,\zeta^{3p-1}X_2,\zeta^{-(3p-2)}X_3)$$
 $$T(X_1,X_2,X_3)=(\eta^{-1}X_3,\eta^{-1}X_1,\eta ^{-1} X_2).$$
 (The notation means, e.g., that $S(X_3)=\zeta^{-(3p-2)}X_3$ and $T(X_2)=\eta^{-1}X_1.$)  Recalling Lemma 1.5, setting  $\gamma=\eta^{3p}$ (a primitive cube root of $1$), one has that 
the polynomial $$f(X_1,X_2,X_3)=X_1^{3p-1}X_2+\gamma X_2^{3p-1}X_3+\gamma^2 X_3^{3p-1}X_1$$ is invariant under the action of $G$.  $f=0$ defines a smooth projective curve $D\subset \mathbb P^2$, containing all the points in Lemma 2.2; the cone $\mathcal D$ over this curve is a hypersurface singularity $\C^3$.  The map $f:\mathbb C^3\rightarrow \C$ is a smoothing of $\mathcal D$, and the Milnor fibre $M$ is simply-connected, with 
Euler characteristic $\chi(M)=1+(3p-1)^3$.  But $G$ acts freely on $\C^3-\{0\}$, and in particular on $M$ and on $\mathcal D-\{0\}$; thus, $f:\C^3/G\rightarrow \C$ gives a smoothing of $\mathcal D/G\equiv V$, with Milnor fibre the free quotient $M/G$.  Since $|G|=mn=\chi(M)$, we conclude that $M/G$ is a $\mathbb Q$HD, with fundamental group $G$.

\begin{theorem}\label{th2.3} Suppose $p\geq 2$, $G$ and $f$ as above.  Then $f:\C^3/G\rightarrow \C$ defines a $\mathbb Q$HD smoothing of a rational surface singularity $V$ of type $\mathcal A^4$, whose Milnor fibre has fundamental group $G$.  The resolution dual graph of $V$ is 

$$
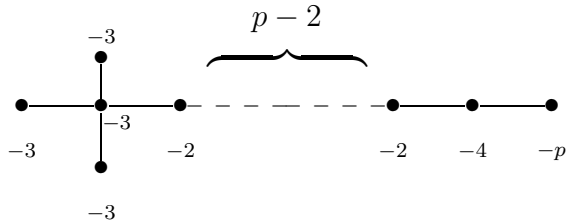
\xymatrix@R=6pt@C=24pt@M=0pt@W=0pt@H=0pt{
\\&&&&p-2\\
&&\overtag{\bt}{-3}{8pt}
&&{\hbox to 0pt{\hss$\overbrace{\hbox to 60pt{}}$\hss}}&\\
&&\lineto[u]\\
&\undertag{\bt}{-3}{6pt}\lineto[r]&
\bt\lineto[u]\lineto[r]_(.2){-3}&\undertag{\bt}{-2}{6pt}\dashto[r]
&\dashto[r]
&\undertag{\bt}{-2}{6pt}\lineto[r]
&\undertag{\bt}{-4}{6pt}\lineto[r]
&\undertag{\bt}{-p}{6pt}\\
&&\lineto[u]\\
&&\lineto[u]\\
&&\undertag{\bt}{-3}{6pt}\lineto[u]&&&}
$$

\vspace{3.5in}
The  set of $4$ points on the central curve is equianharmonic, i.e., the cross ratio is a primitive $6$th root of unity.
\begin{proof}
We outline the proof as described in Example 8.12 of \cite{SSW}; full details are given in the more difficult case of Section $7$.  The blow-up $\mathfrak B\rightarrow \mathcal D$ at the origin is the geometric line bundle over $D$ corresponding to $\mathcal O_D(1)$.   The action of $G$ lifts to $\mathfrak B$; $T^3$ is now a pseudo-reflection, so divide $\mathfrak B$ by it, yielding $\mathfrak B'$, on which $G/<T^3>=\bar{G}$ now acts.  $\mathfrak B'$ still contains a smooth exceptional divisor isomorphic to $D$.  The action of $\bar{G}$ on $\mathfrak B'$ is non-free only at  points of $D$ described in Lemma 2.2.  There are $3$ orbits on which every point has isotropy of order $3$; and three points where the isotropy is generated by $\bar{S}$, and $\bar{T}$ permutes them cyclically.   Consider $4$ points representing these singular orbits; in each case, compute the Jacobian of a generator of the stabilizer, and diagonalize it with respect to the tangent spaces of $\mathfrak B'$ and $D$.  Using (2.7) below, resolve then the $4$ cyclic quotient singularities  of $\mathfrak B'/\bar{G}$, yielding most of the resolution diagram claimed by the Theorem; one lacks only the genus and self-intersection $-d$ of the central curve.  

Riemann-Hurwitz for the quotient of $D$ by $\bar{G}$ yields that the central curve has genus $0$.  The Milnor fibre for the smoothing has as first homology group the abelianization of $G$, hence has order $n$.  As is well-known (e.g., \cite{lw}, 2.4), since the Milnor fibre has no rational homology,  $n$ is the order of a self-isotropic subgroup of the discriminant group.  That implies that the group itself has order $n^2$.    The order of the discriminant group can be computed from the graph ((2.7) below), from which we conclude that $d=3$. 

Finally, one needs to compute the cross-ratio of the $4$ points on $\mathbb P^1$ which are the images of the points on $D$ on which $\bar{G}$ acts with non-trivial isotropy.  

Consider the transformation $U\in GL(3)$ defined by $U(a,b,c)=(a,\gamma b, \gamma^2 c)$  (recall that $\gamma$ is a primitive cube root of $1$).  One checks that $UTU^{-1}$ is $T$ composed with scalar multiplication by $\gamma^2$, which is itself a power or $T$ (either $T^{3p}$ or its inverse).  Thus, $U$ normalizes $G$, and in fact the corresponding $\bar{U}$ centralizes $\bar{G}$.   $\bar{U}$ cyclically permutes the three points $[1,\rho, \rho^2 ]$ on $\mathbb P^2$.  Further, $\bar{U}$ sends the polynomial $f(X_1,X_2,X_3)$ into $\gamma^2f$, so it acts on $D\subset \mathbb P^2$.  Thus, $\bar{U}$ induces an order $3$ automorphism of the quotient $D/\bar{G}=\mathbb P^1$; it fixes the image of $[1,0,0]$, and cyclically permutes the images of the other $3$ branch points.   But a set of $4$ points on $\mathbb P^1$ for which an automorphism of order $3$ fixes one point and cyclically permutes the other $3$  is easily seen to have a cross-ratio equal to a primitive $6$th root of unity.
\end{proof}
\end{theorem}

We remark that $G$ was determined by a specific choice of a primitive $n$th root $\eta$.  On the other hand, the rational surface singularity $V$ of Theorem 2.3 is independent of this choice, since by \cite{p2} and \cite{l}, the singularity must be weighted homogeneous and thus uniquely determined by the location (i.e., cross-ratio) of the four points.  However, a choice of a different root could give a non-conjugate group $G'$, hence a different total space $\C^3/G'$, thus giving a different smoothing for the same singularity.

\begin{theorem}\label{A}  Let $\eta$=exp$(2\pi i/n)$, and denote its inverse by $\eta'$.  Then
\begin{enumerate} 
\item the subgroups $G=G_{\eta}$ and $G'=G_{\eta'}$ of $\text{GL}(3,\C)$ are isomorphic but not conjugate;
\item the singularities $\C^3/G$ and $\C^3/G'$ are not isomorphic, but their links are diffeomorphic;
\item the maps $f:\C^3/G\rightarrow \C$ and $f':\C^3/G'\rightarrow \C$ are distinct one-dimensional smoothing components in the semi-universal deformation of the singularity $V$.
\end{enumerate}
\begin{proof}
The first assertion follows from Theorem 5.5.6 of \cite{wolf}, describing irreducible faithful complex representations $\{\pi_{k,l}\}$ of certain metacyclic groups, and which are fixed point free (off the origin).  His notation is compatible with ours, and his $d=3$ and his $n'=3p$.  Then the group $G_{\eta}$ in our Theorem is the image of the representation $\pi_{m-1,1}$, while $G_{\eta'}$ corresponds to $\pi_{m-1, n'-1}$.  But by Statement (3) of Theorem 5.5.6, the two representations are inequivalent even after composition with an automorphism of the group.  It follows that the subgroups themselves are not conjugate.

A local isomorphism of $\C^3/G$ onto $\C^3/G'$ would lift to an automorphism of $\C^3$ which would intertwine the two group actions.  But this is impossible by the previous remark.  As for the links, one notes from Theorem 5.5.10 of \cite{wolf} that the corresponding orthogonal representations in $\mathbb R^6$ are in fact equivalent; hence, the quotient-spheres $S^5/G$ and $S^5/G'$ are diffeomorphic.

The fact that these one-parameter smoothings are full smoothing components follows from the results in \cite{lw}, and will be discussed in Section 8 below.

\end{proof}
\end{theorem}

\begin{remark} The existence of two $\mu=0$ smoothing components for the singularities in Theorem \ref{A} may be understood as follows.   As in part 1 of \cite{lw}, the first homology group of the link of a rational surface singularity has a natural "discriminant quadratic function", which induces the usual linking pairing.  By Theorem 4.5 of \cite{lw},  a $\mu=0$ smoothing gives rise to a self-isotropic subgroup of the discriminant group; a calculation shows there are two such.  They differ by a diffeomorphism of the link which permutes two of the $4$ points on the central $\mathbb P^1$; but this is not induced by an automorphism of the singularity.  This explains why the total spaces of the smoothings are not isomorphic, but their links are diffeomeorphic.  In (\cite{js}, p. 124), it is claimed that there are $6$ smoothing components, basically because of permutations of the $3$ special points on the $\mathbb P^1$; but the order $3$ automorphism of the singularity acts on each of the two smoothing components.

\end{remark}
\begin{remark}  The metacyclic groups above have $p\geq 2$, but all the calculations make sense for $p=1$.  In this case, the group $G$ is cyclic of order $9$, and one finds the quotient of the cone over a cubic curve in $\mathbb P^2$, whose resolution dual graph is 
$$
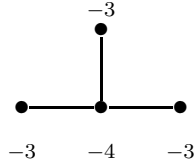
\xymatrix@R=6pt@C=24pt@M=0pt@W=0pt@H=0pt{
\\&&&&\\\
&&\overtag{\bullet}{-3}{8pt}\\
&&\lineto[u]\\
&&\lineto[u]\\
&&\lineto[u]\\
&\undertag{\bullet}{-3}{6pt}\lineto[r]
&\undertag{\bullet}{-4}{6pt}\lineto[u]\lineto[r]&\undertag{\bullet}{-3}{6pt}}
$$
\vspace{.2in}

This is the first term in the $\mathcal W$ series of \cite{SSW}, originally described in (\cite{w3},5.9.2).
\end{remark}
\subsection*{Some familiar facts}
{\bf(2.7)} A two-dimensional cyclic quotient singularity of type $n/q$, where $0<q<n, (q,n)=1$, is the quotient of $\mathbb C^2$ by $(a,b)\mapsto (\eta a, \eta^q b)$, with $\eta$ a primitive \emph{n}th root of $1$.   The minimal resolution is a string of smooth rational curves, the negative of whose self-intersections are the numbers in the continued fraction expansion of $n/q$.  The transform of the image of the coordinate axis $y=0$ intersects transversally the first of these curves.

For a star-shaped resolution dual graph with central curve of self-intersection $-d$ and chains of curves corresponding to continued fraction expansions of $n_i/q_i,\ 1\leq i\leq s$, the order of the discriminant (= absolute value of the determinant of the intersection matrix) is, e.g. via \cite{ne} (1.2),
$$n_1n_2\cdots n_s(d-\sum_{i=1}^{s}q_i/n_i).$$

\section{Class $\mathcal B$ valency 4 smoothings}
\bigskip
Continuing with the previous set-up, we consider representations of dimension $d=4$.  Let $p\geq 2$, and write $m=2p^{2}-2p+1,\ n=8p,$ $n'=2p$, and $r=m-(2p-1).$  Let $\zeta=\text{exp}(2\pi i/m)$ and $\eta=\text{exp}(2\pi i/n)$.   Consider the subgroup $G\subset
GL(4,\C)$ generated by the transformations $S$ and $T$, where 
 $$S(a,b,c,d)=(\zeta a,\zeta^{r}b,\zeta^{-1} c,\zeta^{-r}d)$$ 
 $$T(a,b,c,d)= (\eta b,\eta c,\eta d,\eta a).$$
This is the $4$-dimensional fixed-point free representation $\pi_{1,1}$.
 It is straightforward to verify the following:
 \begin{enumerate}
 \item $T^4$ is scalar multiplication by $\eta^{4}$
 \item The image $\bar{G}$ of $G$ in PGL$(4,\mathbb C)$ is $G/<T^4>$, hence has order $4m$ and is generated by $\bar{S}$ of order $m$ and $\bar{T}$ of order $4$.
 \end{enumerate}
 
 Checking that $r^k-1$ is a unit in $\mathbb Z/(m)$ for $k=1,2,3$, Lemma 1.6 implies
\begin{lemma}
If $1\leq k \leq 3$,  then $S^iT^k$ is conjugate to $T^k.$
\end{lemma}
 
 \begin{lemma}
$\bar{G}$ acts faithfully on $\mathbb P^3$, with non-trivial isotropy only at 
\begin{enumerate}
\item the $4$ points $[1,0,0,0], 
[0,1,0,0], 
[0,0,1,0], [0,0,0,1]$,  each of whose stabilizers is $<\bar{S}>$, and which are cyclically permuted by  $<\bar{T}>$.
\item the $4$ (distinct) $\bar{G}$-orbits of the $4$ points $[1, \rho,\rho^2,  \rho^3]$, where $\rho^4=1$, noting that the stabilizer of each point is  $<\bar{T}>$.
\item  the $\bar{G}$-orbits of the two lines $s[1,0,\mu ,0]+t[0,1,0,\mu ]$, where $\mu=\pm 1$, noting that  $<\bar{T}^2>$ fixes every point on the lines.
\end{enumerate}
\begin{proof}
One checks the action of $\bar{G}$ on $\mathbb P^3$ for orbits with non-trivial isotropy.  By the last Lemma, it suffices to find points fixed by $\bar{S}^i$, $\bar{T}$, and $\bar{T}^2$.  The calculations concerning $\bar{G}$ are now straightforward.

\end{proof}
\end{lemma}
Next, $G$ acts on the 4 coordinate functions on $\C^4$ 
 by $$S(X_1,X_2,X_3,X_4)=(\zeta^{-1} X_1,\zeta^{-r} X_2,\zeta X_3,  \zeta^{r} X_4)$$ 
 $$T(X_1,X_2,X_3,X_4)=(\eta^{-1}X_4,\eta^{-1}X_1,\eta^{-1}X_2,\eta ^{-1} X_3).$$
 Denoting $\gamma=\eta^{2p}$ ($=i$, a primitive $4$th root of $1$), we have:
 
 \begin{lemma} $G$ fixes the polynomial $$f(X_1,X_2,X_3,X_4)=X_1^{2p-1}X_2+\gamma X_2^{2p-1}X_3+\gamma^2X_3^{2p-1}X_4+\gamma^{3} X_4^{2p-1}X_1,$$ and leaves invariant the hypersurface siingularity $$\mathcal Z=\{X_1X_3-X_2X_4=0\}\subset \C^4.$$  Thus $G$ acts on the variety $\mathcal D=\mathcal Z\cap\{f=0\} \subset \C^4$, the  cone over a smooth complete intersection curve $D\subset \mathbb P^3$.
 \begin{proof} The only claim requires checking is the smoothness of the curve in $\mathbb P^3$.  A more difficult but similar argument is given below in Section 6.
 \end{proof}
 \end{lemma}
 
The map  $f:\mathcal Z\rightarrow \C$ gives a smoothing of $\mathcal D$; by \cite{hamm} and \cite{gh}, the Milnor fibre $M$ is simply connected with Euler characteristic which can be checked to equal $mn$.  $G$ acts freely on $\C^4-\{0\}$, and in particular on  $M$ and $\mathcal D-\{0\}$.  Thus, $f:\mathcal Z/G\rightarrow \C$ gives a smoothing of $\mathcal D/G$, with Milnor fibre the free quotient $M/G$.  Since $|G|=mn$, we conclude that $M/G$ is a $\mathbb Q$HD, with fundamental group $G$.

\begin{theorem} Suppose $p\geq 2$, $G$, $\mathcal Z$, and $f$ as above.  Then $f:\mathcal Z/G\rightarrow \C$ defines a $\mathbb Q$HD smoothing of a rational surface singularity $V$ of type $\mathcal B^4$, whose Milnor fibre has fundamental group $G$.  The resolution dual graph of $V$ is 
$$
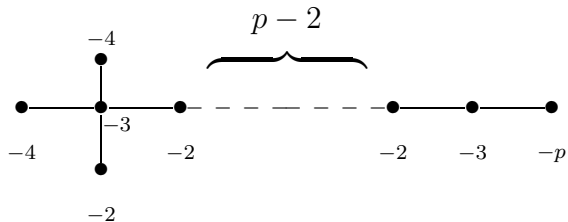
\xymatrix@R=6pt@C=24pt@M=0pt@W=0pt@H=0pt{
\\&&&&p-2\\
&&\overtag{\bt}{-4}{8pt}
&&{\hbox to 0pt{\hss$\overbrace{\hbox to 60pt{}}$\hss}}&\\
&&\lineto[u]\\
&\undertag{\bt}{-4}{6pt}\lineto[r]&
\bt\lineto[u]\lineto[r]_(.2){-3}&\undertag{\bt}{-2}{6pt}\dashto[r]
&\dashto[r]
&\undertag{\bt}{-2}{6pt}\lineto[r]
&\undertag{\bt}{-3}{6pt}\lineto[r]
&\undertag{\bt}{-p}{6pt}\\
&&\lineto[u]\\
&&\lineto[u]\\
&&\undertag{\bt}{-2}{6pt}\lineto[u]&&&}
$$

\vspace{.5in}
The  set of $4$ points on the central curve has cross ratio $2$, as they admit an automorphism of order $2$.
\begin{proof}
Again we outline the proof, as in Theorem \ref{th2.3}.   One resolves the complete intersection singularity $\mathcal D$ by one blow-up, divides out by the action of the pseudo-reflection $T^4$, and locates and describes the fixed points of $\bar{G}$ on this smooth surface using Lemma 3.2.  The fixed points of type $(1)$ yield one $\bar{G}$-orbit, and the isotropy of a point is $<\bar{S}>$, of order $m$.  For type $(2)$, only those two points with $\rho=\pm 1$ are on $\mathcal Z$; the isotropy of each is $<\bar{T}>$, of order $4$.  For type $(3)$, there are only $4$ points on the $2$ lines which lie on $\mathcal Z$, and two of them are of type $(2)$; the remaining points $[1,1,-1,-1]$ and $[1,-1,-1,1]$ are permuted by $\bar{T}$, so give only one $\bar{G}$-orbit, and the isotropy of one is $<\bar{T}^2>$.   From the explicit form of the isotropy, one gets all the information of the resolution dual graph, except for the genus and self-intersection of the central curve.  

The same argument as in the proof of Theorem \ref{th2.3}, using Riemann-Hurwitz and the order of a self-isotropic subgroup of the discriminant group, yields the remaining invariants of the resolution graph.

Finally, one needs to compute the cross-ratio of the $4$ points on $\mathbb P^1$ which are the images of the points on $D$ on which $\bar{G}$ acts with non-trivial isotropy.  

Consider the transformation $U\in GL(4)$ defined by $U(a,b,c,d)=(a,-b,c,-d)$. One checks that $UTU^{-1}=-T=T^{4p+1}$ .  Thus, $U$ normalizes $G$, and in fact the corresponding $\bar{U}$ centralizes $\bar{G}$.   Consider the action of $\bar{U}$ on $4$ points whose orbits yield the $4$ branch points on $\mathbb P^1$.  $\bar{U}$ interchanges $[1,1,1,1]$ and $[1,-1,1,-1]$ and leaves fixed $[1,0,0,0]$ and the $\bar{G}$-orbit of $[1,1,-1,-1]$.  Thus, $\bar{U}$ induces an automorphism of $\mathbb P^1$ which fixes two points and interchanges two others.  It is easily seen that such a set of points has cross-ratio $2$ (or equivalently $1/2$ or $-1$).  
\end{proof}
\end{theorem}

\begin{remark}  The metacyclic groups above have $p\geq 2$, but all the calculations make sense for $p=1$.  In this case, the group $G$ is cyclic of order $8$, and one finds the quotient of the cone over a complete intersection elliptic  curve in $\mathbb P^3$, whose resolution dual graph is 
$$
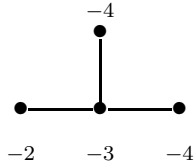
\xymatrix@R=6pt@C=24pt@M=0pt@W=0pt@H=0pt{
\\&&&&\\\
&&\overtag{\bullet}{-4}{8pt}\\
&&\lineto[u]\\
&&\lineto[u]\\
&&\lineto[u]\\
&\undertag{\bullet}{-2}{6pt}\lineto[r]
&\undertag{\bullet}{-3}{6pt}\lineto[u]\lineto[r]&\undertag{\bullet}{-4}{6pt}}
$$
\vspace{.2in}

This is the first term in the $\mathcal N$ series of \cite{SSW}.
\end{remark}

\section{Class $\mathcal C^4$---the group and its representation}
\bigskip

We do a variation of the earlier constructions.  Let $d=6$ and $p\geq 2$, and write $m=p^2-p+1, n=6p, n'=p, $ and $r=m-p+1.$  Constructing transformations $S_k$ and $T_l$ as in $1.4 (1)$ and $1.4(2)$, we have a six-dimensional representation of a semi-direct product $G$ of two cyclic groups, of orders $m$ and $n$.  It is still true that $r$ has order $6$ in $\mathbb Z/(m)$.  However, the representation need not be fixed-point free, either because $m$ and $n$ are not relatively prime when $p\equiv 2$ mod $3$ (their GCD would be $3$), or because the divisors of $n'$ might not include $2$ and $3$.   Not only that, but we are actually interested in a $7$-dimensional \emph{reducible} representation of this group; we take the direct sum of the previous representation with a character of the abelianization of $G$.  

Specifically, write  $\zeta$ = exp$(2\pi i/m)$ and $\eta$ = exp$(2\pi i/n)$.  Define $S,T \in U(7)$ by

$$S(a_1,a_2,\cdots,a_7)=(\zeta a_1,\zeta^{r} a_2,\cdots,\zeta^{r^{5}}a_6,a_7)$$ 
 $$T(a_1,a_2,\cdots,a_7)= (\eta a_2,\cdots, \eta a_{6},\eta a_1,\eta a_7).$$
 Then $S$ and $T$ generate a metacyclic group $G$ satisfying $$S^m=T^n=I, \ TST^{-1}=S^r.$$
 $G$ is a semi-direct product of the cyclic groups $<S>$ and $<T>$, and the abelianization is cyclic of order $n$, generated by the image of $T$.
 
 It is straightforward to verify the following:
 \begin{enumerate}
 \item $T^6$ is scalar multiplication by $\eta^{6}$
 \item The image $\bar{G}$ of $G$ in PGL$(7,\mathbb C)$ is $G/<T^6>$, hence has order $6m$ and is generated by $\bar{S}$ of order $m$ and $\bar{T}$ of order $6$.
 \item For every $i,j,k,l$, one has
 $$(S^{j}T^l)(S^iT^k)(S^{j}T^l)^{-1}=S^{ir^l-j(r^k-1)}T^k.$$

 \end{enumerate}
 
 It is useful to note that \begin{equation*}\tag{4.0.1}S(a_1,a_2,a_3,a_4,a_5,a_6,a_7)=(\zeta a_1,\zeta^{1-p}a_2,\zeta^{-p}a_3,\zeta^{-1}a_4,\zeta^{p-1}a_5,\zeta^{p}a_6,a_7).\end{equation*}
Now, $r-1$,\ $r^3-1$, and $r^5-1$ are units in $\mathbb Z/(m)$,  and so are $r^2-1$ and $r^4-1$ unless $p\equiv 2\  \text{mod}\ 3$.  Using $(3)$ above, we deduce 
 \begin{lemma}\label{le:conj}
\begin{enumerate}
\item If $k$ is congruent to $1,3,\text{or}\ 5\  \text{mod}\  6$, then $S^iT^k$ is conjugate to $T^k$.
\item If $p\not \equiv 2$ mod $3$ and $k$ is congruent to $2$ or $4$ mod $6$, then $S^iT^k$ is conjugate to $T^k$.
\item Suppose $p \equiv 2$ mod $3$.  Then if $k$ is congruent to $2$ or $4$ mod $6$, then $S^iT^k$ is conjugate to $T^k$ if $3|k$ and  $ST^k$ otherwise.
\end{enumerate}
\begin{proof} Calculate that $p=3s+2$ implies that both $r^2-1$ and $r^4-1$ are equivalent modulo $m$ to $3$ times a unit. Also note that $r^3$ is congruent to $-1$ mod $m$.  The Lemma now follows by careful use of (3) above.
\end{proof}
\end{lemma}

We will be interested in the action of $G$ on $\mathbb C^7$ and of $\bar{G}$ on $\mathbb P^6.$ Denote by $e_1,e_2,\cdots,e_7$ the standard basis vectors of $\mathbb C^7$, and by $\bar{e_i}$ the corresponding points in $\mathbb P^6$.
\begin{lemma} Consider the action of $G$ on $\mathbb C^7$.
\begin{enumerate}
\item Every point on the line of multiples of $e_7$ is fixed by the subgroup $<S>$.
\item Any other point fixed where $G$ acts non-freely has last coordinate equal to $0$.
\end{enumerate}
\begin{proof}  If $S^iT^k$ fixes a point whose last coordinate is non-$0$, then $k$ must be $0$.
\end{proof}
\end{lemma}

\begin{lemma}  Consider the action of $\bar{G}$ on $\mathbb P^6$.
\begin{enumerate}
\item $\bar{e_7}$ is fixed by all of $\bar{G}$
\item $<\bar{S}>$ is the stabilizer of each of $\bar{e_1},\cdots,\bar{e_6}$;  $<\bar{T}>$ acts transitively on these points; $<\bar{S}>$ acts freely off $\bar{e_1},\cdots,\bar{e_7}$
\item Any other point where $\bar{G}$ acts non-freely is in the $\bar{G}$-orbit of a point fixed by $\bar{T}^2$ or $\bar{T}^3$ or (if $p\equiv 2\  \text{mod}\ 3$) $\bar{S}\bar{T}^2$.
\end{enumerate}
\begin{proof}  The first two statements are immediate.  To locate orbits on which a group acts non-freely, it suffices to consider fixed points of representatives of each conjugacy class, as described in Lemma 4.1.  Any point fixed by $\bar{T}^k$ is also fixed by either $\bar{T}^2$ or $\bar{T}^3$.   Finally, if $p=3s+2$, one checks that $(ST^4)^2=S^{3s+2}T^8$, which is conjugate via $4.1(3)$ to $ST^8$; so, a fixed point of $\bar{S}\bar{T}^4$ is in the $\bar{G}$-orbit of a fixed point of $\bar{S}\bar{T}^2$.
\end{proof}
\end{lemma}

In fact, we will consider a smooth subvariety $Z\subset \mathbb P^6$ and its affine cone  $\mathcal Z \subset\mathbb C^7$ on which $G$ acts freely off the origin. 
\bigskip
\bigskip

\section{Class $\mathcal C^4$---the variety and the group action}
\bigskip

The linear system of cubics in $\mathbb P ^2$ with base points at $(1,0,0), (0,1,0), (0,0,1)$ gives a rational map $\Phi:\mathbb P^2- -\rightarrow\mathbb P^6$, 
providing a projectively normal embedding of $Z$, the blow-up at these points.   Then $Z\subset \mathbb P^6$ has  $\mathcal O  _Z(1)$ equal to the (very ample) anti-canonical bundle $-K_Z$, and $Z$ has degree $K_Z\cdot K_Z =6$.   Since $H^1(Z,\mathcal O _Z(j))=0$ for all $j$, the
 affine cone $\mathcal Z \subset \mathbb C^7$  has an isolated Cohen-Macaulay three-dimensional singularity at the origin;  it is also canonical Gorenstein, since the canonical bundle of $Z$ is the negative of the hyperplane bundle. 

We find equations for $Z\subset \mathbb P^6$ and for $\mathcal Z$. The rational map $\Phi$ can be given by 
 mapping $[A,B,C]$ to
$$  [X_{1},X_{2},X_3,X_4,X_5,X_6,X_7]\equiv 
[AB^2,A^2B,A^2C,AC^2,BC^2,B^2C,ABC].$$  
Then the homogeneous ideal for  $Z\subset \mathbb P^6$ is generated by $9$ quadrics:

$$X_1X_3-X_2X_7\ \ \ \ \ X_4X_6-X_5X_7\ \ \ \ X_1X_4-X_7^2$$
$$X_2X_4-X_3X_7\ \ \ \ \ X_1X_5-X_6X_7\\ \ \ \ \ X_2X_5-X_7^2$$
$$ X_3X_5- X_4X_7\ \ \ \ \ X_2X_6-X_1X_7\ \ \ \ \ X_3X_6-X_7^2$$

\begin{proposition} Choose $p\geq 2$ and let $G\subset U(7)$ be the group defined in Section 4.  Then $G$ acts on $\mathcal Z \subset \mathbb C^7$ and $\bar{G}$ acts on $Z\subset \mathbb P^6.$
\begin{proof}  Since the matrices in $G$ are unitary, the transpose inverse of such an element is its complex conjugate.  So, the action of $S$ on a coordinate $X_i$ is multiplication by $\zeta^{-r^{i-1}}$, while $T$ sends $X_i$ (for $i\neq 7$) to $\eta^{-1}X_{i-1}$, and $X_7$ to $\eta^{-1}X_7$. Now one can check that $S$ multiplies each of the $9$ defining equations of $\mathcal Z$ by a power of $\zeta$ (it helps to use (4.0.1)), and $T$ sends each equation into $\eta^{-2}$ times another. 
\end{proof}
\end{proposition}
 Let $E_i$, $i=1,\cdots,5,$ denote the line in $\mathbb P^6$ through $\bar{e}_i$ and $\bar{e}_{i+1}$, and by $E_6$ the line through $\bar{e}_6$ and $\bar{e}_1$.   Then the following facts are well-known or easy to verify:
\begin{enumerate}
\item The union of lines $E_i$ forms a cycle of lines on $Z$, and their sum is the divisor $Z\cap\{X_7=0\}.$
\item  Any point on $Z$ for which some coordinate is $0$ lies on this divisor.

\end{enumerate}
   
 \begin{proposition}\label{cor:g} The only points of $Z$ with a coordinate equal $0$ on which $\bar{G}$ acts non-freely are the $6$ points $\bar{e}_1,\cdots,\bar{e}_6$.  For each of these, the stabilizer is $<\bar{S}>$, and these points form one $\bar{G}$-orbit. 
 \begin{proof} The assertions about the $\bar{e}_i$ are clear.  One must show that any other point $\bar{P}$ on, e.g.,  $E_1$ cannot be fixed by a non-trivial $\bar{S}^i\bar{T}^j$.  But $\bar{P}$ has non-0 coordinates exactly in the first two entries.   If it is fixed by  $\bar{S}^i\bar{T}^j$, then $j=0$.   But, applying $\bar{S}^i$ with $i>0$ multiplies the entries by two different powers of $\zeta$, hence changes the point. \end{proof}
 \end{proposition}
 
  Note $\bar{G}$ acts on $Z-\cup E_i$.  To study  fixed points of $\bar{G}$ on this set, it suffices to look at points on $Z$ within the standard affine open in $\mathbb P^6$ where $X_1\neq 0$.   Considering coordinates $X_i/X_1$, $i=2,\cdots ,7,$ we find using the $9$ defining equations that $U=Z\cap\{X_1\neq 0\}$ is here an affine plane with coordinates  $x=X_2/X_1$ and $y=X_6/X_1.$  So,  points of $U$ may be written as $$(a,b)=[1,a,a^2b,a^2b^2,ab^2,b,ab];$$
 In these $x,y$ coordinates, $\bar{e}_1$ is the origin, $E_6$ is given by $x=0$, and  $E_1$ is given by $y=0$.  The action of $\bar{G}$ on  a point $(a,b)$ in $U-\{xy=0\}$ is given by
  \begin{equation}\label{eq:g}\bar{S}(a,b)=(\zeta^{-p}a,\zeta^{p-1}b)
  \end{equation}
  \begin{equation}\label{eq:h} \bar{T}(a,b)=(ab,1/a).
 \end{equation}

 It is clear that the cyclic group $<\bar{S}>$ acts freely off the origin.   Lemma \ref{le:conj} restricts the transformations to be considered.

\begin{lemma}\label{lem:x}  Consider the action of $\bar{G}$ on $U-\{(a,b)|ab=0\}$.
\begin{enumerate}
\item The only fixed point of $\bar{T}$ is $\bar{Q}_1=(1,1).$
\item The fixed points of $\bar{T}^2$ are the $3$ points $(\rho,\rho),$ where $\rho^3=1.$
\item If $p=3s+2$, then  the fixed points of $\bar{S}\bar{T}^2$  are the $3$ points $(\rho \zeta^{-(s+1)},\rho \zeta^{2s+1})$, where $\rho^3=1$.

\item The fixed points of $\bar{T}^3$  are the $4$ points $(\pm 1,\pm 1).$ 
\end{enumerate}
\begin{proof} These are simple calculations, using the explicit form of $\bar{S}$ and $\bar{T}$.
\end{proof}
\end{lemma}

\begin{proposition}\label{pr:2} The points of $Z-\cup E_i$ one which $G$ acts non-freely are the orbits in $U$ of:
\begin{enumerate}
\item $\bar{Q}_1=(1,1)$, with stabilizer $<\bar{T}>$ of order $6$,
 \item If $p\not \equiv 2$ mod $3$, the point $\bar{Q}_2=(\omega, \omega)$, where $\omega$= exp$(2\pi i/3)$, with stabilizer $<\bar{T}^2>$ of order $3$,
 \item If $p=3s+2$, the point $\bar{Q}_3=(\omega \zeta^{-(s+1)},\omega \zeta^{2s+1})$, with stabilizer $<\bar{S}\bar{T}^2>$ of order $3$.
 \item $\bar{Q}_4=(-1,-1),$ with stabilizer $<\bar{T}^3>$ of order 2.
 \end{enumerate}
 
 \begin{proof}
 The orbits of the fixed points are exactly the orbits of the points in Lemma \ref{lem:x}; one needs to enumerate the distinct orbits, and find the stabilizer of a representative fixed point in each orbit.  The stabilizer of $\bar{Q}_1$ is easily checked to be as asserted.
     The $3$ points in Lemma \ref{lem:x}(2) are  $\bar{Q}_1$, $\bar{Q}_2$, and $\bar{T}(\bar{Q}_2)$.  A calculation shows that if $p\not \equiv 2$ mod $3$, then the stabilizer of $\bar{Q}_2$ is $<\bar{T}^2>$ , and it is not in the $\bar{G}$-orbit of $\bar{Q}_1$.   On the other hand, if $p=3s+2$, then one can write $m=3m'$ and $\bar{Q}_2=\bar{S}^{m'}(\bar{Q}_1)$, so these points are already accounted for in the orbit of $\bar{Q}_1$.  But then the $3$ points in Lemma \ref{lem:x}(3) are $\bar{Q}_3$, $\bar{S}^{m'}(\bar{Q}_3)$, and $\bar{S}^{2m'}(\bar{Q}_3)$; another calculation shows the stabilizer of $\bar{Q}_3$ is as claimed, and the point is not in one of the preceding orbits.  Finally, the $4$ points in Lemma \ref{lem:x}(4) are
 $\bar{Q}_1$, $\bar{Q}_4$, $\bar{T}(\bar{Q}_4)$, and $\bar{T}^{2}(\bar{Q}_4)$, and the calculation of the stabilizer of $\bar{Q}_4$ is straightforward.    
 \end{proof}
 \end{proposition}
 
  We  summarize and complete the above results in the

  \begin{theorem}\label{th:1} 
 \begin{enumerate}
 \item The metacyclic group $G$ acts freely on $\mathcal Z -\{0\}$, hence 
  $\mathcal Z /G$ has an isolated $\mathbb Q$-Gorenstein (in particular, Cohen-Macaulay) singularity.
  \item The induced action of $\bar{G}$ on $Z$ is faithful, and has $4$ orbits of fixed points 
    \end{enumerate}
  \begin{proof} We show $G$ acts freely on $\mathcal Z -\{0\}$.  We already know the fixed points for the induced action on $\bar{Z}$.  If $S^iT^j$  has an eigenvector in $\mathcal Z$ with eigenvalue $1$, then consideration of the last coordinate shows that either it is $0$, or $j=0$.   In the first case, by Proposition 5.2 it suffices to consider $e_1$.  It is an eigenvector  for $S^iT^{6j},\ i<m$, with eigenvalue $\zeta^{i} \eta^{6j}$; but a simple check shows this cannot equal $1$ unless the transformation is the identity.  The case $j=0$ again yields only the $e_i$ as possible fixed points, and $e_7$ is not on $\mathcal Z$.  
  
  The other assertions follow from Proposition \ref{pr:2}.

  \end{proof}
  \end{theorem}
\bigskip
\bigskip

\section{Class $\mathcal C^4$---the $G$-invariant polynomial}
\bigskip

We continue to consider for each $p\geq2$ the subgroup $G\subset U(7)$ from the previous sections. 
Write $\gamma=\eta^{p} = \text{exp}(2\pi i/6)$, a primitive sixth root of $1$, and consider the polynomial $$f=X_1^{p-1}X_2+\gamma X_2^{p-1}X_3+\gamma^2X_3^{p-1}X_4+\gamma^3X_4^{p-1}X_5+\gamma^4X_5^{p-1}X_6+\gamma^5X_6^{p-1}X_1.$$

 
 \begin{theorem} $f$ is $G$-invariant, and defines a smooth curve $D$ on $Z$.  The induced action of $\bar{G}$ on $D$ is faithful, and the quotient is a $\mathbb P^1$, with $D\rightarrow \mathbb P^1$ branched over $4$ points.
 \begin{proof}
 A direct computation shows that $f$ is $G$-invariant.  The (possibly non-reduced) curve $D\subset Z$ which it defines contains all the points of $Z$ where $\bar{G}$ has isotropy, as enumerated in Propositions \ref{cor:g} and \ref{pr:2} .
 
We consider $D$ on the open subset $U$ of $Z$ given by $\{X_1\neq 0\}$.   As in the discussion earlier, one has coordinates $x$ and $y$ for $Z$ for which  $\bar{e}_1$ is the origin, $E_6$ is given by $x=0$, and  $E_1$ is given by $y=0$.  $D$ is the plane curve defined by the vanishing of
  $$f(x,y)=
  x +\gamma x^{p+1}y
 +\gamma^2x^{2p}y^{p+1}
  +\gamma^3x^{2p-1}y^{2p}
  +\gamma^4x^{p-1}y^{2p-1}
  +\gamma^5y^{p-1}.$$
$D$ intersects the two axes only at the origin, and it is smooth there; it intersects $E_6$ with multiplicity $p-1$, and $E_1$ with multiplicity $1$.  

We show $D$ is reduced and irreducible.   Any irreducible component of $D$ must intersect the ample divisor $\{X_7=0\}$; the above description shows it can intersect only at the $\bar{e_i}$ ($i\neq 7$), and it is smooth there; so, $D$ is reduced.  Let $C$ be the irreducible component of $D$ containing $\bar{e_1}$.  Since $\bar{T}$ permutes the $\bar{e}_i$,  the set of curves $\{\bar{T}^j(C),\ j=0, \cdots,5\}$ consists of all irreducible components of $D$; so, there are at most $6$ of them.

We prove in Lemma \ref{D} below that $D$ contains all the points at which $\bar{G}$ acts non-freely, and is smooth there.  In particular, if $D$ had a singular point, then it would have at least $|\bar{G}|=6(p^2-p+1)$ of them; let us assume this.  If $\tilde{D}$ is the normalization of $D$, then the standard exact sequence yields
 $$\chi(\mathcal O_{\tilde{D}})=\chi(\mathcal O_D)+\sum \delta_Q,$$ the last sum adding up the $\delta$-invariants of the singular points of $D$.  Since  $\chi(\mathcal O_D)=-D\cdot (D+K)/2$ and $D\equiv pH$, the previous formula yields the inequality
 $$h^0(\mathcal O_{\tilde{D}})-h^1(\mathcal O_{\tilde{D}})\geq -3p(p-1)+6(p^2-p+1).$$
 But we have observed that $h^0(\mathcal O_{\tilde{D}})\leq 6$; then the inequality $h^1(\mathcal O_{\tilde{D}})\geq 0$ yields the contradiction $0\leq 3p-3p^2.$

Consider finally the action of $\bar{G}$ on the smooth curve $D$.   A simple topological count (or Hurwitz's formula) shows that  the genera of $D$ and the quotient $D/\bar{G}$ are related by
$$(2g(D)-2)/ |\bar{G}|\  -\sum_{i=1}^{s}(1-1/n_i)=2g(D/\bar{G})-2,$$
where there are $s$ orbits of fixed points, and $n_i$ is the order of the  stabilizer of a representative point.   But $2g(D)-2=6p(p-1)$ and $s=4$, with the $n_i$ equal $2,3,6,$ and $m$; plugging in yields that $D/\bar{G}$ is a projective line.
\end{proof}

\begin{lemma}\label{D}  The points of $Z$ at which $\bar{G}$ acts non-freely are all smooth points of $D$.  Further, $D$ has the following tangent lines for certain points in $U$:
\begin{enumerate}
\item At $\bar{Q}_1=(1,1)$: $(x-1)+\gamma(y-1)=0.$ 
\item At $\bar{Q}_2=(\omega, \omega)=(\gamma ^2, \gamma ^2)$ : $(x-\omega)+\gamma(y-\omega)=0.$
\item At $\bar{Q}_3=(\omega \zeta^{2s+1},\omega \zeta^{-s})$, when $p=3s+2$: $(x-\omega \zeta^{2s+1})+\gamma \zeta^{-p}(y-\omega \zeta^{-s})=0.$
\end{enumerate}
\begin{proof}  Since $\bar{G}$ acts on $D$, it suffices to check the first assertion at a representative of each of the $4$ orbits where the action is not free.  $\bar{e}_1$ is obviously a point of $D$, and smoothness there has already been noted.  The verification that $\bar{Q}_i $ for $i=1,2,3$ are points of the curve, with $f_y/f_x$ as claimed at each point, is a very lengthy computation using high-school algebra.  It is easy to check that $f(-1,-1)=0$ and $f_x(-1,-1)=2-2p-2p(-1)^p+\gamma(-1)^p(2p-2)\neq 0$; it will not be necessary to know the equation of the tangent line at $\bar{Q}_4$.
\end{proof}
\end{lemma}

 \end{theorem}

  \section{A $\mathbb Q$HD smoothing for $\mathcal C^4$}
  
  Consider as before the surface $Z\subset \mathbb P^6$, its affine cone $\mathcal Z \subset \mathbb C^7$, and for each $p\geq 2$ the group $G=G_p$, and the $G$-invariant homogeneous polynomial $f$ of degree $p$.  The map $f:\mathcal Z\rightarrow \mathbb C$ has  special fibre $\mathcal D$,  the affine cone over the smooth curve $D\subset \mathbb P^6$ given by $Z\cap\{f=0\}$.  As $K_D=\mathcal O_D(p-1)$, we see that $\mathcal D$ has an isolated Gorenstein singularity, and the map $f$ expresses $\mathcal Z$ as the total space of a smoothing.  The Milnor fibre $M$ for this smoothing is any fibre of $f:\mathcal Z -f^{-1}(0)\rightarrow \mathbb C-\{0\}.$   
    
    \begin{proposition} The Milnor fibre $M$ for this smoothing of $\mathcal D$ is simply connected, with Euler characteristic $6p(p^2-p+1)$.
    \begin{proof} Since $\mathcal D$ is Gorenstein, the Euler characteristic of its Milnor fibre $M$ is given by the familiar formula computed from a resolution
    $$1+\mu=12p_g+c_1^2+c_2.$$
 The genus of $D$ is $3p(p-1)+1$, and  one easily checks that $c_1^2=-6p^3,\  c_2=-6p(p-1), \ p_g=p^3$, from which the desired formula follows.
    
    For simple connectivity of the Milnor fibre, the usual Lefschetz theory shows that it suffices to prove simple connectivity of the link of the smoothing, i.e., of $L=\mathcal Z \cap S^{13}\subset \mathbb C^7$.  But $L$ is homeomorphic to the $S^1$-bundle on $Z$ corresponding to the hyperplane line bundle; so, by the exact sequence of homotopy groups and the simple connectivity of $Z$, we see $\pi_1(L)$ is abelian.   It suffices therefore to prove the vanishing of the first homology group of $L$; we use the Gysin sequence, part of which is
    $$H_2(Z)\rightarrow H_0(Z)\rightarrow H_1(L)\rightarrow H_1(Z)=0.$$
    The first map is cupping with the first Chern class of the hyperplane bundle, so its cokernel is a cyclic group whose order is the degree of imprimitivity of $c_1$ in $H^2(Z)$, or equivalently in Pic($Z$).  But the hyperplane class is primitive (e.g., because its self-intersection of $6$ is square-free). This proves the claim.  
    \end{proof}
    \end{proposition}
    
    \begin{remark}  A Milnor fibre $M=f^{-1}(t)$ of the smoothing can be constructed as the complement of a curve on a projective surface.  Consider in $\mathbb P^7$ the variety $Z'$ defined by the equations of $Z$ plus a new equation $T^p-f=0$ (where $T$ is a new homogeneous coordinate). Note $Z'\cap \{T=0\}=D$, and $Z'\rightarrow Z$ is a cyclic covering ramified over $D$.   The Milnor fibre is easily seen to be $Z'-D$, a cyclic $p$-fold unramified covering of $Z-D$, and the Euler characteristic could be computed directly in this way.
    \end{remark}

 Since $G$ acts freely on the cones $\mathcal Z$ and $\mathcal D$, and $f$ is $G$-invariant, we have an induced map $f:\mathcal Z /G\rightarrow \mathbb C$, a smoothing of the surface singularity $\mathcal D/G$.  
 
 \begin{proposition} $G$ acts freely on the Milnor fibre $M$ above, and the quotient $M/G$ is the Milnor fibre of the smoothing $f$ of $\mathcal D/G$.  The $4$-manifold $M/G$ is a rational homology disk.
 \begin{proof} Since $G$ acts freely on $\mathcal Z-\{0\}$ and $f$ is $G$-invariant, we have that $G$ acts freely on $M=f^{-1}(t), t\neq 0$.  So, the Euler characteristics satisfy $\chi(M)=|G|\chi(M/G)$.  Via Proposition $7.1$, we conclude that $\chi(M/G)=1$.  But $\pi_1(M/G)\simeq G$, since $M$ is simply-connected, and so $H^i(M/G;\mathbb Q)=0$ for $i>0$.  Thus $M/G$ is a rational homology disk.
 \end{proof}
 \end{proposition}
 
 We are now ready to state the main theorem.
  
   \begin{theorem}\label{C}  Suppose $p\geq 2$, with $G$, $\mathcal Z$, and $f$ as in Sections 4, 5, and 6 respectively.  Then $f:\mathcal Z /G \rightarrow \mathbb C$ defines a $\mathbb Q$HD smoothing of a rational surface singularity $V$ of type $\mathcal C^4$, whose Milnor fibre has fundamental group $G$.  The resolution dual graph of $V$ is 
 
$$
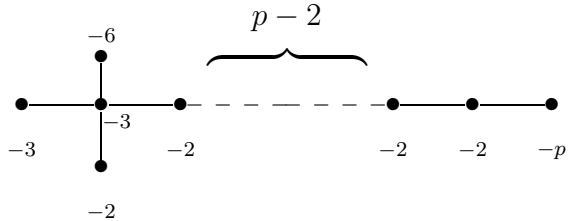
\xymatrix@R=6pt@C=24pt@M=0pt@W=0pt@H=0pt{
\\&&&&p-2\\
&&\overtag{\bt}{-6}{8pt}
&&{\hbox to 0pt{\hss$\overbrace{\hbox to 60pt{}}$\hss}}&\\
&&\lineto[u]\\
&\undertag{\bt}{-3}{6pt}\lineto[r]&
\bt\lineto[u]\lineto[r]_(.2){-3}&\undertag{\bt}{-2}{6pt}\dashto[r]
&\dashto[r]
&\undertag{\bt}{-2}{6pt}\lineto[r]
&\undertag{\bt}{-2}{6pt}\lineto[r]
&\undertag{\bt}{-p}{6pt}\\
&&\lineto[u]\\
&&\lineto[u]\\
&&\undertag{\bt}{-2}{6pt}\lineto[u]&&&}
$$
\vspace{1in}

\begin{proof}
Given the Propositions above, it remains only to find the resolution dual graph of the singularity $\mathcal D/G$.  

First, consider the blowing-up $\mathfrak B\rightarrow \mathcal Z$ at the origin, yielding the geometric line bundle over $Z$ associated to $\mathcal O_Z(-1)$; the action of $G$ lifts to this space.  The non-free orbits occur at the points of $Z\subset \mathfrak B$  described in Proposition 5.4; as every such orbit contains a point with first coordinate non-$0$, it suffices to consider the open set $\mathcal U\subset \mathfrak B$ obtained by inverting the coordinate $X_1$.  As in the proof of Lemma $5.3$,  $\mathcal U$ is a $3$-dimensional affine space with coordinates $x=X_2/X_1,\ y=X_6/X_1, z=X_1$, and  the exceptional divisor $Z$ is given by $z=0$.  $G$ acts on $\mathcal U-\{xy=0\}$.
 The actions of $S$ and $T$ on $\mathcal U$ are given by 
 $$S(a,b,c)=(\zeta^{-p}a, \zeta^{p-1}b, \zeta  c),$$
 $$T(a,b,c)=(ab,1/a,\eta ac),$$
 hence $$T^6 (a,b,c)=(a,b,\eta^6 c).$$
Since $T^6$ is a pseudo-reflection (of order $p$), we divide $\mathcal U$ by it, getting a cyclic map $\mathcal U\rightarrow \mathcal U'=\mathcal U/<T^6>$ to another affine space, branched $p$ times over the last coordinate.   $\mathcal U'$ has coordinates $x,y,z'=z^p$, with $z'=0$ giving the image of $Z$ (the exceptional divisor), and $f(x,y)=0$  (as in (6.1)) the proper transform of $D$.  $\bar{G}=G/<T^6>$ acts on $\mathcal U'$ via
  $$\bar{S}(a,b,c')=(\zeta^{-p}a, \zeta^{p-1}b, \zeta^{p} c'),$$
 $$\bar{T}(a,b,c')=(ab,1/a,\eta^pa^pc').$$
 The relevant points on $\mathcal U'$ to consider and their $\bar{G}$-stabilizers are:
 \begin{enumerate}
 \item $(0,0,0)$, with stabilizer $<\bar{S}>$ of order $m$.
 \item $(1,1,0)$, with stabilizer $<\bar{T}>$ of order $6$.
 \item If $p \not \equiv 2$ mod $3$, $(\omega,\omega,0)$, with stabilizer $<\bar{T}^2>$ of order $3$.
 \item If $p=3s+2$, $(\omega \zeta^{-(s+1)},\omega \zeta^{2s+1},0)$, with stabilizer $<\bar{S}\bar{T}^2>$ of order $3$.
 \item $(-1,-1,0)$, with stabilizer $<\bar{T}^3>$, of order $2$.
 \end{enumerate}
 
These data would allow a complete description of the four three-dimensional cyclic quotient singularities on $\mathcal U'/\bar{G}$, or equivalently on the partial resolution $\mathfrak B/G.$ 
 
But we are interested in the quotient singularities on a partial resolution of $\mathcal D/G$.  For this, we consider in $\mathcal U'$ the smooth two-dimensional hypersurface $f(x,y)=0$, and divide by $\bar{G}$. 
At each of the special points above, one considers the Jacobian matrix of a generator of the stabilizer, and examines its action on the tangent space to $f(x,y)=0$ (given by the formulas of Lemma 6.2), now paying attention to the exceptional divisor $z'=0$ corresponding to $D$.  Following $(2.7)$, we describe each cyclic quotient singularity by a fraction $n/q$, counting out from the "central curve" which is the image of $D$.
 
 At $(0,0,0)$, the tangent plane is given by $x=0$.  For a point in this plane, the action of $\bar{S}$ is $(b,c')\mapsto(\zeta^{p-1}b,\zeta^pc')$.   Note that $\zeta^p=(\zeta^{p-1})^{m-p+1}$; so, resolving and starting from the image of $D$, one has an $m/m-p+1$ cyclic quotient singularity.
 
 At $(1,1,0)$, the tangent plane is $(x-1)+\gamma(y-1)=0$.  From the Jacobian matrix of $\bar{T}$, one finds the  cyclic quotient is of type $6/1$.  The cases with stabilizer of order $3$ are done similarly, again using the previous computation of the tangent plane; one has a cyclic quotient of type $3/1$.  Finally, for the point $(-1,-1,0)$, the stabilizer $\bar{T}^3$ has order $2$, hence to first order acts as minus the identity, which gives an $A_1$ on $D/\bar{G}$.
 
 These calculations give most of the information on the resolution dual graph of $\mathcal D/G$.  That the central curve is rational was proved above in Theorem $6.1$.  There remains only the question of the self-intersection $-d$ of the central curve.  As in $(2.7)$, the discriminant equals $$2\cdot 3\cdot 6\cdot m\cdot \{d-(1/2)-(1/3)-(1/6)-(m-p+1/m)\}=6^2\{(d-3)(p^2-p+1)+p^2\}.$$  On the other hand, as the singularity has a $\mu=0$ smoothing, there is a self-isotropic subgroup of the discriminant group which is isomorphic to the first homology of the Milnor fibre (e.g., \cite{lw}, 4.5).  As noted above, that first homology group is the abelianization of $G$, hence has order $n=6p$.  Thus, the order of the discriminant group is $n^2=(6p)^2$, so $d=3$.

\end{proof}
\end{theorem}

\begin{remark}  The metacyclic groups above have $p\geq 2$, but again all the calculations make sense for $p=1$.  In this case, the group $G$ is cyclic of order $6$, and one finds the quotient of the cone over a hyperplane section of $Z\subset \mathbb P^6$ (an elliptic curve of degree $6$), whose resolution dual graph is 
$$
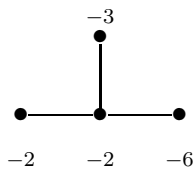
\xymatrix@R=6pt@C=24pt@M=0pt@W=0pt@H=0pt{
\\&&&&\\\
&&\overtag{\bullet}{-3}{8pt}\\
&&\lineto[u]\\
&&\lineto[u]\\
&&\lineto[u]\\
&\undertag{\bullet}{-2}{6pt}\lineto[r]
&\undertag{\bullet}{-2}{6pt}\lineto[u]\lineto[r]&\undertag{\bullet}{-6}{6pt}}
$$
\vspace{.2in}

This is the first term in the $\mathcal M$ series of \cite{SSW}.
\end{remark}

\section{$\mathbb Q$HD smoothing components}

 An irreducible component of the base space of the semi-universal deformation of a normal surface singularity $(V,0)$ is called a \emph{smoothing component} if the generic fibre over such a component is smooth.  For a rational surface singularity, every component is a smoothing component, but their dimensions vary.  We are interested in components containing a $\Q$HD smoothing.
\begin{theorem}  A $\mathbb Q$HD smoothing component of a rational surface singularity has dimension $$h^1(S_X)+\sum_{i=1}^{r} (d_i-3),$$
where $(X,E)\rightarrow (V,0)$ is the minimal good resolution; $E=\sum_{i=1}^{r}E_i$;  $S_X$, the sheaf of logarithmic vector fields, is the kernel of the natural surjection $\Theta_X\rightarrow \oplus N_{E_i}$; and $E_i \cdot E_i=-d_i$.
\begin{proof} In general, suppose that $\pi:(\mathcal V,0)\rightarrow (\mathbb C,0)$ is a smoothing of a normal surface singularity, and let $\Theta_{\mathcal V /\mathbb C}$ denote the relative derivations.  Consider the length $\beta$ of the cokernel of the natural map
$$\Theta_{\mathcal V /\mathbb C} \otimes \mathcal O_V\rightarrow \Theta_V.$$
It was conjectured in \cite{w3}(4.2) that $\beta$ is the dimension of the corresponding smoothing component; this was then proved by Greuel-Looijenga \cite{gl} .   Next, assuming the smoothing could be globalized in an appropriate sense, Theorem 3.13 of \cite{w3} proved a formula for $\beta$ in terms of a resolution of $V$ and the Euler characteristic of the Milnor fibre.  The globalization condition (called $(G)$ in that theorem) was later shown by Looijenga (\cite{lo}, Appendix) to be satisfied in all cases.  
Applying the formula of Theorem 3.13(c) of \cite{w3} to the case at hand, one has $h^1(\mathcal O_X)=0$ and $\chi_T(E)=1+r$.  Thus, for  a $\mathbb Q$HD smoothing, the corresponding component has dimension $$\beta=h^1(\Theta_X) -2r.$$
Note that $S_X$ is the dual of the perhaps more familiar $\Omega^1_X(\text{log} E)$ (see \cite{w1},(2.2)).  Since $h^0(N_{E_i})=0$, and $h^1(N_{E_i})=d_i-1$ if $E_i\cdot E_i =-d_i$, one obtains the stated formula. 
\end{proof}
\end{theorem}

For a rational surface singularity, $h^1(S_X)$ is the space of first-order deformations of $X$ preserving all the exceptional curves; it is the tangent space to the smooth functor of equisingular deformations of $V$ (\cite{w1},(5.16)).  Further, for every exceptional cycle $Z$ there is a surjection $S_X\rightarrow \Theta_Z$, so that the induced surjective map $H^1(S_X)\rightarrow H^1(\Theta_Z)$ is an isomorphism for all $Z$ sufficiently ``large".

On the other hand, a resolution graph is ``taut" if there is a unique analytic singularity with this graph.  H.  Laufer classifies all of these in \cite{l}.  By standard theorems, this is equivalent to $H^1(\Theta_Z)=0$ for all exceptional cycles $Z$ on a resolution of any singularity with this graph, so that a taut singularity will have $H^1(S_X)=0$.  

\begin{corollary}  A $\mathbb Q$HD smoothing component of a weighted homogeneous surface singularity has dimension 1.
\begin{proof}  It is easily checked that all graphs of types $\mathcal G, \mathcal W, \mathcal N, \mathcal M$ of \cite{SSW} are taut, and satisfy $\sum_{i=1}^{r} (d_i-3)=1$.  This settles the assertion in those cases.  For graphs of types  $\mathcal A, \mathcal B, \mathcal C$ with one star, is follows from the allowed moves as described in Figure $11$ of \cite{SSW} that $\sum_{i=1}^{r} (d_i-3)=1$ for all final graphs of valency $3$, and $\sum_{i=1}^{r} (d_i-3)=0$ for final graphs of valency 4.   Again, in the valency $3$ case, Laufer's result shows all these graphs are taut.

    For the valency $4$ graphs, H. Laufer has shown (\cite{l}, (4.1)) that the analytic structure is uniquely determined by the cross-ratio of the $4$ points on the central curve, i.e., by the analytic structure of the reduced exceptional cycle $E$.  Thus, $h^1(\Theta_Z)=h^1(S_X)=1$ in all cases.  Combining with our remark above that $\sum_{i=1}^{r} (d_i-3)=0$, the assertion about the dimension of the component is proved. 
\end{proof}
\end{corollary}

\begin{corollary}  The smoothings in Theorems 2.3, 3.4, and 7.4 give smooth, one-dimensional smoothing components in the base space of the versal deformation of the corresponding singularities.
\begin{proof}  Since the smoothing components have dimension $1$, one needs only to check that the Kodaira-Spencer map on the tangent spaces is injective, i.e., the deformations of the surface singularities are non-trivial to first order.  We do this in two steps.  Write the singularity $V$ as a quotient $Y/G$.  First, the deformation of $Y$ is non-trivial to first order.  This follows from the easily proved general result:

\begin{lemma} Let $R=\mathbb C[[x_1,\ldots,x_n]]/(f_i)$ be a local domain of embedding dimension $n$, $h\in m_R$ a non-zero divisor.  Then the map $\mathbb C[[t]]\rightarrow R$ sending $t$ to $h$ is flat, giving a one-parameter deformation of $R/(h)$.  If $h\in m_R^2$, then the induced first-order deformation of $R/(h)$is non-trivial, i.e., the Kodaira-Spencer map is non-$0$.
\end{lemma}

Thus, one has a non-trivial first order deformation of the surface singularity $(Y,0)$. We claim next that the corresponding deformation of $(Y/G,0)$ is non-trivial to first-order.  But this follows from the injectivity of $T^1_{Y/G}\subset (T^1_Y)^G,$ which can be proved using the injectivity of $H^1(Y/G-\{0\},\Theta_{Y/G})\subset H^1(Y-\{0\},\Theta_Y)^G$.
\end{proof}
\end{corollary}

\begin{remark} The same argument gives smooth one-dimensional smoothing components for other examples of quotient $\Q$HD smoothings, as in \cite{SSW}, (8.2).
\end{remark}

\begin{theorem}  Consider an $H$-shaped graph of type $\mathcal A, \mathcal B$, or $\mathcal C$, i.e.
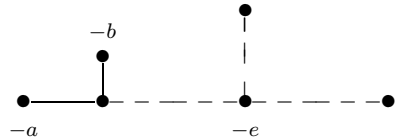
$$
\xymatrix@R=12pt@C=24pt@M=0pt@W=0pt@H=0pt{
&&&&{\bullet}\dashto[dd]\\
&&\overtag{\bullet}{-b}{10pt}&&\dashto[u]\\
&\undertag{\bullet}{-a}{0pt}\lineto[r]&
{\bullet}\lineto[u]\dashto[r]&\dashto[r]&
\undertag{\bullet}{-e}{0pt}\dashto[r]&\dashto[r]&{\bullet}\\
&&&&
}
$$
Thus, $a$ and $b$ are any two of the $3$ integers in one of the three triples.  Suppose that the graph is taut, e.g., $e\geq 3$.  Then there is no singularity with this resolution graph having a $\mathbb Q$HD smoothing.
\begin{proof}  The inductive definition of a graph of type $\mathcal A, \mathcal B$, or $\mathcal C$ consists of two types of blow-up moves on a preliminary graph, as described in Figure $11$ of \cite{SSW}.  One obtains the desired graph at any stage by replacing the (unique) $-1$ node by respectively a  $-4,\ -3, \ -2$.  It is easy to see that the second type of move leaves  $\sum_{i=1}^{r} (d_i-3)$ the same, while the first type lowers it by $1$.  Starting with the basic preliminary graph,  this sum equals $1$; to construct an $H$-shaped graph also requires a series of type $(b)$ moves, and then one type $(a)$ move, followed perhaps by more type $(b)$ moves.  The net result is that an $H$-shaped graph of any of these types will have $\sum_{i=1}^{r} (d_i-3)=0$.  But a taut graph would have $h^1(S_X)=0$, so Proposition 8.1 implies one has a $\Q$HD smoothing component of dimension $0$, an impossibility.

         Table IV of \cite{l} lists the taut $H$-shaped graphs (all of which are rational save one case which could not arise here).  If $e\geq 3$ in the $H$-shape, the graph above is taut, of type $(L_1)+(J_1)+(R_2)$.  
\end{proof}
\end{theorem}

There are many taut graphs with $e=2$ which are of type $\mathcal A, \mathcal B, \mathcal C$; so, none of these could have a $\Q$HD smoothing.   

\begin{corollary} Consider the graphs $Y_n$ from Figure 4 of \cite{gs}:
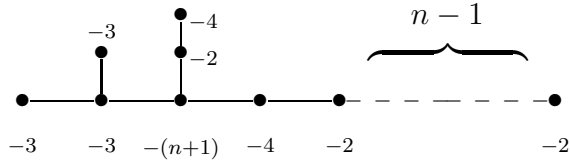
$$
\xymatrix@R=6pt@C=24pt@M=0pt@W=0pt@H=0pt{
\\&&&\righttag{\bullet}{-4}{6pt}
&&&n-1\\
&&\overtag{\bullet}{-3}{8pt}&\righttag{\bullet}{-2}{6pt}\lineto[u]
&&&{\hbox to 0pt{\hss$\overbrace{\hbox to 60pt{}}$\hss}}&\\
&&\lineto[u]&\lineto[u]\\
&\undertag{\bullet}{-3}{6pt}\lineto[r]&
\undertag{\bullet}{-3}{6pt}\lineto[u]\lineto[r]&\undertag{\bullet}{-(n+1)}{6pt}\lineto[u]\lineto[r]&\undertag{\bullet}{-4}{6pt}\lineto[r]&\undertag{\bullet}{-2}{6pt}\dashto[r]&\dashto[r]
&\undertag{\bullet}{-2}{6pt}\\
&\\
&\\
&\\
&}
$$
They are taut for all $n\geq 1$, hence the singularities do not have a $\Q$HD smoothing.
\begin{proof}  The case  $n\geq 2$ is covered by the condition above $e\geq 3$.  For $n=1$, the graph is taut, as it is of type $(L_1)+(J_2)+(R_4)$, using the contraction $(C_3)$ to end up with the desired graph.
\end{proof}
\end{corollary}

\begin{remark}  By the method of ``symplectic caps", Gay and Stipsicz proved (\cite{gs}, Proposition 4.2) that for $n\geq 7$, the $3$-manifolds constructed from $Y_n$ do not bound a $\Q$HD, hence \emph{a fortiori} the corresponding singularities could not have a $\Q$HD smoothing.
\end{remark}
\begin{remark} The following graph (of type $\mathcal A$) is rational but not taut, so the above argument does not rule out that with some analytic structure, there exists a $\Q$HD smoothing:
$$
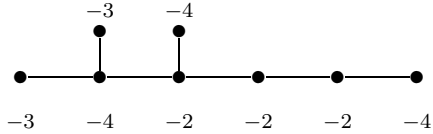
\xymatrix@R=6pt@C=24pt@M=0pt@W=0pt@H=0pt{
\\&&&&\\
&&\overtag{\bullet}{-3}{8pt}&\overtag{\bullet}{-4}{8pt}
&&&\\
&&\lineto[u]&\lineto[u]\\
&\undertag{\bullet}{-3}{6pt}\lineto[r]&\undertag
{\bullet}{-4}{6pt}\lineto[u]\lineto[r]&\undertag{\bullet}{-2}{6pt}\lineto[u]\lineto[r]&\undertag{\bullet}{-2}{6pt}\lineto[r]
&\undertag{\bullet}{-2}{6pt}\lineto[r]
&\undertag{\bullet}{-4}{6pt}\\
&&\\
&&\\
&&&&&\\
&\\
&}
$$
We note that this graph passes the ``$\bar{\mu}=0$" test as well (cf. \cite{siegen}, \cite{smu}), so cannot be eliminated on those grounds either.
\end{remark}
Nonetheless, we make the following
\begin{conjecture}The only complex surface singularities admitting a $\mathbb Q$HD smoothing are the (known) weighted homogeneous examples.
\end{conjecture}

The question of existence of examples other than those on our original list was asked even before the organization of the examples in \cite{SSW}, see for instance the remark on the bottom of page $505$ of \cite{ds}.

\end{document}